%% file: ModalLogicOfForcing.tex
\documentclass[12pt]{amsart}
\usepackage{diagrams}
\diagramstyle[tight,centredisplay%
,dpi=600%              for more modern laser-printers
,noPostScript%
,heads=LaTeX%
]%
\input ArticleMacros
\renewcommand{\implies}{\to}
\renewcommand{\iff}{\leftrightarrow}
\newcommand{\Force}{\mathrm{Force}}
\def\S#1{\hbox{$\mathsf{S#1}$}}
\newcommand{\axiomf}[1]{{\rm #1}}
\newcommand{\theoryf}[1]{\hbox{$\mathsf{#1}$}}
\begin{document}
\author{Joel David Hamkins}
\address{The College of Staten Island of The City University of New York, Mathematics, 2800 Victory Boulevard, Staten Island, NY 10314, USA
\& The Graduate Center of The City University of New York, Ph.D.~Program in Mathematics, 365 Fifth Avenue, New York, NY 10016, USA}
\email{jhamkins@gc.cuny.edu, http://jdh.hamkins.org}

\author{Benedikt L\"owe}
\address{Institute for Logic, Language and Computation, Universiteit van Amsterdam,
Plantage Muidergracht 24, 1018 TV Amsterdam, The Netherlands} \email{bloewe@science.uva.nl}

\bottomnote In addition to partial support from PSC-CUNY grants and other CUNY support, the first author was a \textit{Mercator-Gastprofessor} at the
Westf\"alische Wilhelms-Universit\"at M\"unster in Summer 2004, when this collaboration began, and was partially supported by \textit{NWO
Bezoekersbeurs} \textsf{B 62-612} at Universiteit van Amsterdam in Summer 2005, when it came to fruition. The second author was partially supported by
\textit{NWO Reisbeurs} \textsf{R 62-605} during his visits to New York and Los Angeles in January and February 2005. The authors would like to thank
Nick Bezhanishvili (Amsterdam), Dick de Jongh (Amsterdam), Marcus Kracht (Los Angeles CA), and Clemens Kupke (Amsterdam) for sharing their knowledge of
modal logic.

\Title The modal logic of forcing

\Abstract A set theoretical assertion $\psi$ is {\df forceable} or {\df possible}, written $\possible\psi$, if $\psi$ holds in some forcing extension,
and {\df necessary}, written $\necessary\psi$, if $\psi$ holds in all forcing extensions. In this forcing interpretation of modal logic, we establish
that if \ZFC\ is consistent, then the \ZFC-provable principles of forcing are exactly those in the modal theory $\S{4.2}$.
% When parameters are allowed, some of the possibilities have large cardinal consistency strength.

\Section Introduction

What are the most general principles in set theory relating forceability and truth? We are interested in how the set theoretical method of forcing
affects the first order theory of a model of set theory. As with Solovay's celebrated analysis of provability, both this question and its answer are
naturally formulated with modal logic.\footnote{\cite{Solovay1976:ProvabilityModalLogic}; for a survey of the result and the subsequent development of
the field of provability logic, see also \cite{JaparidzeDeJongh1998:LogicOfProvability}.} We aim to do for forceability what Solovay did for
provability.

Forcing was introduced by Paul Cohen in 1962 in order to prove the independence of the Axiom of Choice \AC\ and the Continuum Hypothesis \CH\ from the
other axioms of set theory. In an explosion of applications, set theorists subsequently used it to construct an enormous variety of models of set
theory and prove many other independence results. With forcing one builds an extension of any model $V$ of set theory, in an algebraic manner akin to a
field extension, by adjoining a new ideal object $G$, a $V$-generic filter over a partial order $\P$ in the ground model $V$, while preserving \ZFC.
The resulting forcing extension $V[G]$ is closely related to the ground model $V$, but may exhibit different set theoretical truths in a way that can
often be carefully controlled. The method has become a fundamental tool in set theory.

Because the ground model $V$ has some access via names and the forcing relation to the objects and truths of the forcing extension $V[G]$, there are
clear affinities between forcing and modal logic. (One might even imagine the vast collection of all models of set theory, related by forcing, as an
enormous Kripke model.) Accordingly, we define that a statement of set theory $\varphi$ is {\df forceable} or {\df possible} if $\varphi$ holds in some
forcing extension, and $\varphi$ is {\df necessary} if it holds in all forcing extensions. The modal notation $\possible\varphi$ and
$\necessary\varphi$ expresses, respectively, that $\varphi$ is possible or necessary. This forcing interpretation of modal logic was introduced by the
first author in \cite{Hamkins2003:MaximalityPrinciple} in connection with the Maximality Principle, a new forcing axiom, with related work in
\cite{Leibman2004:Dissertation} and \cite{HamkinsWoodin:NMPccc}. An alternative but related connection between modal logic and forcing was explored by
Fitting and Smullyan in \cite{FittingSmullyan1996:SetTheoryAndCH}, and Blass \cite{Blass1990:InfinitaryCombinatoricsAndModalLogic} provides an
interpretation of modal logic in set theory that is not directly related to forcing.

These modal operators, of course, are eliminable in the language of set theory, because their meaning can be expressed in the usual language of set
theory by means of the forcing relation or Boolean values. For example, $\possible\varphi$ simply means that there is some partial order $\P$ and
condition $p\in\P$ such that $p\forces_\P\varphi$, and $\necessary\varphi$ means that for all partial orders $\P$ and $p\in\P$ we have
$p\forces_\P\varphi$. In this way, one can interpret $\possible\varphi$ and $\necessary\varphi$ in any model of set theory.\footnote{By formalizing
forcing in \ZFC, rather than the metatheory, one can sensibly force over any model of \ZFC, without needing it to be countable or transitive and
regardless of the metatheoretical objects, such as generic filters, which may or may not exist in a larger universe. In this syntactic account, one
considers what is forced by various conditions without ever building the forcing extension as a structure. A semantic account of forcing over an
arbitrary $M\satisfies\ZFC$ is provided by the quotient $M^\P/U$ of the Boolean valued universe $M^\P$ by any ultrafilter $U$, with no need for $U$ to
be $M$-generic (even $U\in M$ works fine!); since $M$ maps elementarily into the ground model of $M^\P/U$, one has the forcing extension as an actual
structure.} In \ZFC\ we may freely use a mixed language of set theory with the modal operators $\possible$ and $\necessary$, understood with the
forcing interpretation. While the modal operators $\possible$ and $\necessary$ are eliminable, we nevertheless retain them, because our goal is to
discover which modal principles forcing must obey.

For example, it is easy to see that $\necessary\varphi\implies\varphi$ is a valid principle of forcing, because if $\varphi$ is true in all forcing
extensions, then it is true, as the universe is a (trivial) forcing extension of itself. Similarly, $\neg\possible\varphi\iff\necessary\neg\varphi$ is
valid for forcing because a statement $\varphi$ is not forceable if and only if all forcing extensions satisfy $\neg\varphi$. The principle
$\necessary\varphi\implies\necessary\necessary\varphi$ is valid because if $\varphi$ holds in all forcing extensions, then so does $\necessary\varphi$,
since any forcing extension of a forcing extension is a forcing extension. The reader may easily verify that
$\necessary(\varphi\implies\psi)\implies(\necessary\varphi\implies\necessary\psi)$ is valid. The principle
$\possible\necessary\varphi\implies\necessary\possible\varphi$ is valid for forcing, because if $\varphi$ is necessary in $V^\P$ and $V^\Q$ is an
arbitrary extension, then $\varphi$ is true in the product extension $V^{\P\cross\Q}$, as this extends $V^\P$; consequently, $\varphi$ is forceable
over every such extension $V^\Q$. These modal assertions axiomatize the modal theory known as \S{4.2}, and a bit of formalization will help us express
what we have observed. A {\df modal assertion} is a formula of propositional modal logic, such as $(\necessary q_0\implies q_0)$, expressed with
propositional variables $q_i$, Boolean connectives $\wedge$, $\vee$, $\neg$, $\implies$, $\iff$ and the modal operators $\possible, \necessary$.

\Definition Main Definition. A modal assertion $\varphi(q_0,\ldots,q_n)$ is a {\df valid principle of forcing} if for all sentences $\psi_i$ in the
language of set theory, $\varphi(\psi_0,\ldots,\psi_n)$ holds under the forcing interpretation of $\possible$ and $\necessary$.

\noindent More specifically, $\varphi(q_0,\ldots,q_n)$ is a {\df \ZFC-provable principle of forcing} if \ZFC\ proves all such substitution instances
$\varphi(\psi_0,\ldots,\psi_n)$. This naturally generalizes to larger theories with the notion of a $T$-provable principle of forcing. For any model
$W\satisfies\ZFC$, the modal assertion $\varphi(q_0,\ldots,q_n)$ is a {\df valid principle of forcing in $W$} if all substitution instances
$\varphi(\psi_0,\ldots,\psi_n)$ are true in $W$. So $\varphi$ is a valid principle of forcing if it is valid in the universe $V$ of all sets (this is
expressed as a scheme), and the \ZFC-provable principles of forcing are those provable in \ZFC, as schemes, to be valid.

\Question Main Question. What are the valid principles of forcing? \label{MainQuestion}

For the \ZFC-provable principles of forcing, our Main Theorem \ref{MainTheorem} provides the answer. Meanwhile, a few paragraphs ago, we essentially
observed the following.

\Theorem. Every assertion in the modal theory \S{4.2} is a \ZFC-provable principle of forcing.\label{S4.2IsValid}

\Proof. The modal theory \S{4.2} is obtained from the axioms
$$\begin{array}{rl}
 \axiomf{K} & \necessary(\varphi\implies \psi)\implies(\necessary\varphi\implies\necessary\psi)\\
 \axiomf{Dual} & \neg\possible \varphi\iff \necessary\neg\varphi\\
 \axiomf{S} & \necessary\varphi\implies\varphi\\
 \axiomf{4} & \necessary\varphi\implies \necessary\necessary\varphi\\
 \axiomf{.2} & \possible\necessary\varphi\implies\necessary\possible\varphi\\
\end{array}$$
by closing under modus ponens and necessitation (from $\varphi$, deduce $\necessary\varphi$). We observed earlier that each axiom is a valid principle
of forcing. The \ZFC-provable principles of forcing are clearly closed under modus ponens, and they are closed under necessitation because if
$\varphi(\psi_0,\ldots,\psi_n)$ holds in all models of set theory, then so does $\necessary\varphi(\psi_0,\ldots,\psi_n)$.\QED

Let us quickly show that none of the modal theories most commonly considered beyond \S{4.2} are \ZFC-provable principles of forcing. This follows from
our main theorem, but we find it informative to have explicit failing instances of the principles; they are what pointed to \S{4.2} as the right
choice.
$$\begin{array}{rl}
 \axiomf{5} & \possible\necessary\varphi\implies\varphi\\
 \axiomf{M} & \necessary\possible\varphi\implies\possible\necessary\varphi\\
 \axiomf{W5} & \possible\necessary\varphi\implies(\varphi\implies\necessary\varphi)\\
 \axiomf{.3} & \possible\varphi\wedge\possible\psi\implies(\possible(\varphi\wedge\possible\psi)
             \vee\possible(\varphi\wedge\psi)\vee\possible(\psi\wedge\possible\varphi))\\
 \axiomf{Dm} & \necessary(\necessary(\varphi\implies\necessary\varphi)\implies\varphi)
               \implies(\possible\necessary\varphi\implies\varphi)\\
 \axiomf{Grz} & \necessary(\necessary(\varphi\implies\necessary\varphi)\implies\varphi)\implies\varphi\\
 \axiomf{\hbox{L\"ob}} & \necessary(\necessary\varphi\implies\varphi)\implies\necessary\varphi\\
 \axiomf{H} & \varphi\implies\necessary(\possible\varphi\implies\varphi)\\
\end{array}$$

\Observation. If\/ \ZFC\ is consistent, then none of the above assertions is a \ZFC-provable principle of
forcing.\label{Observation.NotValidBeyondS4.2}

\Proof. (\axiomf5) Let $\varphi$ be the assertion ``$\omega_1^\mathrm{L}$ is countable,'' where $\mathrm{L}$ denotes G\"odel's constructible universe.
Since the class $\mathrm{L}$ is absolute to all forcing extensions, the expression $\omega_1^\mathrm{L}$ refers to the same ordinal in all forcing
extensions. Since forcing can collapse $\omega_1$ (see \cite[14.3]{Jech:SetTheory3rdEdition}), and it can never be uncollapsed,
$\possible\necessary\varphi$ holds in any model of \ZFC. But $\varphi$ is false in $\mathrm{L}$, so 5 is not valid for forcing over $\mathrm{L}$.

(\axiomf{M}) This is the McKinsey axiom, also known as .1. The Continuum Hypothesis \CH\ is forceable over any model of set theory (see
\cite[ex.~15.15]{Jech:SetTheory3rdEdition}), so $\necessary\possible\CH$ holds. But $\neg\CH$ is also forceable
(\cite[14.32]{Jech:SetTheory3rdEdition}), so $\possible\necessary\CH$ fails. Thus, \axiomf{M} is invalid for forcing in every model of \ZFC.

(\axiomf{W5}) Let $\varphi$ be the assertion ``$\omega_1^\mathrm{L}<\omega_1$ or $\CH$,'' which is true in $\mathrm{L}$, but not necessary there, since
one can force $\neg\CH$ without collapsing $\omega_1^\mathrm{L}$. Since one could collapse $\omega_1^\mathrm{L}$, however, $\varphi$ is forceably
necessary in any model of set theory. Thus, $\possible\necessary\varphi$ holds in $\mathrm{L}$, but not $\varphi\implies\necessary\varphi$, and so
$\mathrm{W5}$ fails as a principle of forcing in $\mathrm{L}$.

(\axiomf{.3}) Let $\varphi$ assert ``$\omega_1^\mathrm{L}$ is collapsed, but not $\omega_2^\mathrm{L}$'' and $\psi$ assert ``$\omega_2^\mathrm{L}$ is
collapsed, but not $\omega_1^\mathrm{L}$.'' These are each forceable over $\mathrm{L}$, since $\omega_1^\mathrm{L}$ and $\omega_2^\mathrm{L}$ can be
collapsed independently. Yet, whenever one of them holds, the other becomes unforceable. So the antecedent of this instance of \axiomf{.3} is true in
$\mathrm{L}$, while the conclusion is false, so .3 is not valid in $\mathrm{L}$.

(\axiomf{Dm}) This axiom is also known as \axiomf{Dum}. Let $\varphi$ assert ``$\omega_1^\mathrm{L}<\omega_1$ or $\CH$,'' as in case ($\mathrm{W5}$).
We observed that $\varphi$ is possibly necessary, but false in $\mathrm{L}$, so $\mathrm{L}\models\neg(\possible\necessary\varphi\implies\varphi)$,
falsifying the conclusion of this instance of $\mathrm{Dm}$. For the antecedent, suppose a forcing extension satisfies $\mathrm{L}^\P\models
\omega_1^\mathrm{L} = \omega_1$. Since $\CH$ is forceable over $\mathrm{L}^\P$ without collapsing $\omega_1^\mathrm{L}$, we conclude
$\mathrm{L}^\P\models \neg\necessary(\varphi\implies\necessary\varphi)$. It follows vacuously that
$\mathrm{L}^\P\satisfies\necessary(\varphi\implies\necessary\varphi)\implies\varphi$. Alternatively, if $\mathrm{L}^\P\models \omega_1^\mathrm{L} <
\omega_1$, then $\varphi$ is true there, so $\mathrm{L}^\P\satisfies\necessary(\varphi\implies\necessary\varphi)\implies\varphi$. Thus, every forcing
extension satisfies $\necessary(\varphi\implies\necessary\varphi)\implies\varphi$, and so
$\mathrm{L}\models\necessary(\necessary(\varphi\implies\necessary\varphi)\implies\varphi)$, fulfilling the antecedent of \axiomf{Dm} in $\mathrm{L}$.
So \axiomf{Dm} is not valid for forcing over $\mathrm{L}$.

(\axiomf{Grz}) This is the Grzegorczyk axiom. Since \axiomf{Grz} is stronger than \axiomf{Dm}, it is invalid in $\mathrm{L}$ by the previous case. We
will show, however, that \axiomf{Grz} is invalid in every model of set theory. Let $\varphi$ be the \CH. Since $\neg\CH$ is forceable over any model of
set theory, $\necessary\varphi$ is always false, and so $\varphi\implies\necessary\varphi$ is equivalent to $\neg\varphi$. So this instance of
\axiomf{Grz} reduces to $\necessary(\necessary\neg\CH\implies\CH)\implies\CH$. Since $\CH$ is always forceable, $\necessary\neg\CH$ is false, making
the first implication, and the hypothesis of this instance of $\mathrm{Grz}$, vacuously true. So this axiom will fail whenever \CH\ fails. Similar
reasoning shows that if $\varphi$ is $\neg\CH$, then the axiom fails if $\CH$ holds. So \axiomf{Grz} fails as a principle of forcing in every model of
set theory.

(L\"ob) This axiom expresses the content of L\"ob's theorem in provability logic, where one interprets $\necessary\varphi$ as ``$\varphi$ is
provable.'' Under the forcing interpretation, however, it is invalid. If $\varphi$ is \CH, then $\necessary\varphi$ is always false, so
$\necessary(\necessary\varphi\implies\varphi)$ is always true. So in any model of set theory, the hypothesis of this instance of the L\"ob axiom is
true, while the conclusion is false. So L\"ob is never a valid principle of forcing.

(H) Let $\varphi$ be the \CH. As we have mentioned, this is necessarily possible, so $\possible\varphi$ always holds. In this case, therefore,
$\necessary(\possible\varphi\implies\varphi)$ is equivalent to $\necessary\varphi$, which is false. So $\mathrm{H}$ fails in any model where \CH\ is
true. By using $\neg\CH$ in the other models, we conclude that $\mathrm{H}$ fails as a principle of forcing in every model of set theory.\QED

The corresponding modal theories, listed below with their implication diagram, are obtained by closing the axioms under modus ponens, substitution and
necessitation. This list is not exhaustive, as there are continuum many modal theories above \S{4.2} that are not listed. We refer the reader to
\cite{ChagrovZakharyaschev1997:ModalLogic} and \cite{BlackburnDeRijkeVenema2001:ModalLogic} for excellent developments of modal logic, including the
analysis of these and many other theories.

$$\quad\hbox{\tiny\vbox{\hbox{Some Common Modal Theories}\hbox{ }\hbox{ }\hbox{$\begin{array}{rcl}
 \S5    &=& \S4+5  \\
 \S{4W5}&=& \S4+\axiomf{W5}\\
 \S{4.3}&=& \S4+\axiomf{.3}\\
 \S{4.2.1}&=& \S4+\axiomf{.2}+\axiomf{M}\\
 \S{4.2}&=& \S4+\axiomf{.2}\\
 \S{4.1}&=& \S4+\axiomf{M}\\
 \S4     &=& \theoryf{K4}+\axiomf{S}\\
 \theoryf{Dm.2}&=& \S{4.2}+\axiomf{Dm}\\
 \theoryf{Dm}&=& \S4+\axiomf{Dm}\\
 \theoryf{Grz}&=& \theoryf{K}+\axiomf{Grz}=\S4+\axiomf{Grz}\\
 \theoryf{GL}&=& \theoryf{K4}+\axiomf{\hbox{L\"ob}}\\
 \theoryf{K4H}&=& \theoryf{K4}+\axiomf{H}\\   % note that this is the *new* H axiom
 \theoryf{K4} &=& \theoryf{K}+4\\
 \theoryf{K} &=& \axiomf{K}+\axiomf{Dual}\\
\end{array}$}}}\!\!\hbox{\tiny\begin{diagram}[height=2em]
            &           &           &           &  \S5      &           &             &            &        \\
            &           &           &           &  \dTo     &           &             &            &        \\
            &           &           &           &  \S{4W5}  &           &             &            &        \\
            &           &           &           &  \dTo     &  \rdTo    &             &            &        \\
            &           & \S{4.2.1} &           & \S{4.3}   &           &\mathsf{Dm.2}&            & \mathsf{Grz} \\
            &           & \dTo      & \rdTo     &  \dTo     &  \ldTo    &   \dTo      & \ldTo      &            & \mathsf{K4H}       \\
 \mathsf{GL}&           & \S{4.1}   &           & \S{4.2}   &           &\mathsf{Dm}  &            & \ldTo(5,5) &        \\
            & \rdTo(4,4)&           & \rdTo     &  \dTo     &  \ldTo    &             &            &  \\
            &           &           &           &  \S4      &           &             &            &        \\
            &           &           &           & \dTo      &           &             &            &       \\
            &           &           &           &\mathsf{K4}&           &             &            &        \\
            &           &           &           &  \dTo     &           &             &            &        \\
            &           &           &           &\mathsf{K} &           &             &            &        \\
\end{diagram}}$$

\Corollary. If\/ \ZFC\ is consistent, then none of the modal theories \S5, \S{4W5}, \S{4.3}, \S{4.2.1}, \S{4.1}, \theoryf{Dm.2}, \theoryf{Dm},
\theoryf{K4H}, \theoryf{GL} or \theoryf{Grz} are \ZFC-provable principles of forcing, and all are invalid in $\mathrm{L}$. The modal theories \S{4.1},
\S{4.2.1}, \theoryf{K4H}, \theoryf{GL} and \theoryf{Grz} are invalid as principles of forcing in every model of set theory.

So if the \ZFC-provable principles of forcing constitute any previously known modal theory, then the best remaining candidate is \S{4.2}.

\bigskip
This article is intended primarily for two audiences: set theorists interested in the fundamental principles of forcing and modal logicians interested
in the application of their subject to set theory. While we felt it necessary in the arguments to assume a basic familiarity with forcing, we do
provide references to specific elementary forcing results in the standard set theory graduate textbooks where this might be helpful. We were able to
provide in our arguments a complete account of the necessary concepts from modal logic. Sections \ref{Section.TheMainQuestion} and
\ref{Section.TheJankovFineFormula} cover the proof of our Main Theorem \ref{MainTheorem}, which answers our Main Question \ref{MainQuestion} above.
Section \ref{Section.TheMainQuestion} contains a complete proof of the main theorem written primarily with the set theoretical reader in mind, and
Section \ref{Section.TheJankovFineFormula} emphasizes certain aspects of the proof for the modal logicians. After this, we apply our technique
fruitfully in Sections \ref{Section.ForcingOverFixedModel}--\ref{Section.RestrictingToAClass} to various other instances of Main Question 2, by
restricting the focus to a given model of set theory, by investigating the role of parameters in the valid principles of forcing, and by restricting
attention to a natural class of forcing notions, such as those with the countable chain condition.

\Section The Main Theorem \label{Section.TheMainQuestion}

Our main theorem provides an answer to Question \ref{MainQuestion}.

\Theorem Main Theorem. If\/ \ZFC\ is consistent, then the \ZFC-provable principles of forcing are exactly those in the modal theory
\S{4.2}.\label{MainTheorem}

The rest of this section is devoted to the proof of this theorem, beginning with the key concepts of buttons and switches. A {\df switch} is a
statement $\varphi$ of set theory such that both $\varphi$ and $\neg\varphi$ are necessarily possible, so that by forcing $\varphi$ can be switched on
or off at will. For example, the \CH\ is a switch, because you can ensure either \CH\ or $\neg\CH$ by forcing over any model of set theory. In
contrast, a {\df button} is a statement that is (necessarily) possibly necessary. The button is {\df pushed} when it is necessary, and otherwise {\df
unpushed}. The idea is that you can always push a button by making it necessary, but having done so, you cannot unpush it again. The assertion
``$\omega_1^\mathrm{L}$ is countable'' is a button because it can be forced over any model of set theory and once it becomes true it remains true in
all further extensions. This button is unpushed in $\mathrm{L}$. The reader is invited to check that a statement is possibly necessary if and only if
it is necessarily possibly necessary; at bottom, this amounts to the \S{4.2} validity of
$\possible\necessary\varphi\iff\necessary\possible\necessary\varphi$. Thus, a button remains a button in every forcing extension. Although it may seem
at first that buttons and switches are very special sorts of statements, we invite the reader to check that in fact every statement in set theory is
either a button, the negation of a button, or a switch (and these types are disjoint).

A collection of buttons $b_n$ and switches $s_m$ is {\df independent} in a model if first, all the buttons are unpushed in the model and second,
necessarily (that is, in any forcing extension), any of the buttons can be pushed and any of the switches can be switched without affecting the value
of any of the other buttons or switches. In other words, the collection of buttons and switches is independent in $W$ if the buttons are unpushed in
$W$, but in any forcing extension $W^\P$, whatever the pattern of buttons and switches is in $W^\P$, any button can be turned on by forcing to some
$W^{\P*\Qdot}$ without affecting the value of any of the other buttons or switches and any switch can be turned either on or off by forcing to some
$W^{\P*\Rdot}$ without affecting the value of any of the other buttons or switches. It follows, of course, that any finite pattern of buttons and
switches being on or off is possible by forcing over $W$.

We note that the counterexample substitution instances showing the forcing invalidity of the modal assertions in Observation
\ref{Observation.NotValidBeyondS4.2} were each Boolean combinations of independent buttons and switches.

\SubLemma. If\/ $\mathrm{V=L}$, then there is an independent collection of infinitely many buttons and infinitely many
switches.\label{Lemma.ButtonsExist}

\Proof. Let the button $b_n$ be the assertion ``$\omega_n^\mathrm{L}$ is not a cardinal,'' and let the switch $s_m$ be the assertion ``the \GCH\ holds
at $\aleph_{\omega+m}$.'' If $\mathrm{V=L}$, then it is clear that none of the buttons is true, but in any model of set theory, the button $b_n$ can be
made true by collapsing $\omega_n^\mathrm{L}$, without affecting the truth of any other  button or the properties of the \GCH\ above $\aleph_\omega$
(see \cite[15.21]{Jech:SetTheory3rdEdition}). Once the button $b_n$ becomes true, it is clearly necessary, because the cardinal $\omega_n^\mathrm{L}$
will remain collapsed in any further forcing extension. The switches $s_m$ are clearly switches, because with forcing one can arrange the values of the
continuum function at $\aleph_{\omega+m}$ at will by forcing that adds no new bounded sequences below $\aleph_\omega$ (see \cite[15.18 \&
related]{Jech:SetTheory3rdEdition}). Thus, the switches can be set to any desired pattern without affecting any of the buttons.\QED

A {\df lattice} is a partial order such that any two nodes $a$ and $b$ have a greatest lower bound or {\df meet}, denoted $a\wedge b$ and a least upper
bound or {\df join}, denoted $a\vee b$. It follows that every finite set $A$ has a meet $\bigwedge A$ and join $\bigvee A$, and that a finite lattice
has a least and a greatest element.

\SubLemma. If\/ $F$ is a finite lattice and $W$ is a model of set theory with a sufficiently large independent family of buttons $b_i$, then to each
node $w\in F$ we may assign an assertion $p_w$, a Boolean combination of the buttons, such that $W$ satisfies:
\begin{enumerate}
 \item In any forcing extension, exactly one of the $p_w$ is true. And $W\satisfies p_{w_0}$, where $w_0$ is the minimal node of $F$.
 \item In any forcing extension satisfying $p_w$, the statement $p_v$ is forceable if and only if $w\leq v$ in $F$.
\end{enumerate}\label{Lemma.LatticeLabels}

\Proof. Let us associate a button $b_u$ with each node $u\in F$. For any $A\of F$, let $b_A=(\bigwedge_{u\in A}\necessary b_u)\wedge(\bigwedge_{u\notin
A}\neg\necessary b_u)$ be the sentence asserting that exactly the buttons in $A$ are pushed and no others.  Let $p_w=\bigvee\set{b_A\st w=\bigvee A}$
be the sentence asserting that the pattern of buttons that have been pushed corresponds to a set $A$ with least upper bound $w$ in $F$. Since every
forcing extension must have some pattern $A$ of buttons $b_u$ pressed and every such $A$ has a least upper bound in $F$, it is clear that $p_w$ will be
true in the model if and only if $w$ is the least upper bound of $A$. Thus, in any forcing extension exactly one of the $p_w$ is true. And $W\satisfies
p_{w_0}$, as all buttons are unpushed in $W$.

For the second claim, suppose that $W[G]$ is a forcing extension where $p_w$ holds. Let $A=\set{u\in F\st W[G]\satisfies b_u}$ be the set of buttons
that are pushed in $W[G]$. Since $p_w$ is true, it must be that $w$ is the join of $A$ in $F$. If $w\leq v$ in $F$, then by pushing the button $b_v$
and no others, we arrive at a forcing extension $W[G][H]$ with buttons pushed in $A\union\singleton{v}$. Since this has join $v$, this means that $p_v$
is true in $W[G][H]$, and hence $p_v$ is forceable in $W[G]$, as desired. Conversely, suppose that $p_v$ is forceable in $W[G]$. Thus, there is some
further extension $W[G][H]$ satisfying $p_v$. This extension exhibits some pattern of buttons $\set{b_u\st u\in B}$, where the join of $B$ is $v$.
Since $W[G][H]$ is a forcing extension of $W[G]$, the buttons in $A$ remain pushed in $W[G][H]$, and so $A\of B$. Thus, $v$ must be at least as large
as the join of $A$, which is $w$, so $w\leq v$ in $F$.\QED

A {\df pre-lattice} is obtained from a lattice by replacing each node with a cluster of one or more equivalent nodes, all related by $\leq$.
Equivalently, it is a partial pre-order $\leq$ (a reflexive and transitive relation) on a set $F$, such that the quotient of $F$ by the equivalence
relation $u\equiv v\iff u\leq v\leq u$ is a lattice under the induced quotient relation $\leq$.

\SubLemma. If\/ $F$ is a finite pre-lattice and $\set{b_i,s_j}_{i,j}$ is a sufficiently large finite independent family of buttons and switches in a
model of set theory $W$, then to each $w\in F$ we may assign an assertion $p_w$, a Boolean combination of the buttons and switches, such that $W$
satisfies:\label{Lemma.PreLatticeLabels}%
\begin{enumerate}
 \item In any forcing extension, exactly one of the $p_w$ is true. And $W\satisfies p_{w_0}$ for any desired node $w_0$ in the minimal cluster of $F$.
 \item In any forcing extension satisfying $p_w$, the statement $p_v$ is forceable if and only if $w\leq v$ in $F$.
\end{enumerate}

\Proof. The idea is to use the buttons as in Lemma \ref{Lemma.LatticeLabels} to determine which cluster is intended in the quotient lattice, and then
use the switches to determine which node is intended within this cluster. Let $[u]$ denote the equivalence class of $u$ in the quotient lattice
$F/{\equiv}$, and let $p_{[u]}$ be the label assigned to $[u]$ in Lemma \ref{Lemma.LatticeLabels}. Thus, $p_{[u]}$ is the disjunction of various
complete patterns of buttons having supremum $[u]$. Suppose that the largest cluster of $F$ has $k$ nodes, and $k\leq 2^n$. For each subset $A\of
\set{0,\ldots,n-1}=n$, let $s_A=(\bigwedge_{i\in A}s_i)\wedge(\bigwedge_{i\notin A}\neg s_i)$ assert that the pattern of switches is specified by $A$.
Since every pattern of switches is possible by forcing over $W$, every $s_A$ is necessarily possible, and in any forcing extension of $W$, exactly one
$s_A$ holds. For each cluster $[u]$, assign to every $w\in[u]$ a nonempty set $\vec A_w=\set{A_0^w,\ldots,A_{j_w}^w}$ of subsets of $n$ in such a way
that the various $\vec A_w$ for $w\in[u]$ partition the subsets $A\of n$. We may assign the pattern $A$ of switches that happen to hold in $W$ to any
desired node $w_0$ in the minimal cluster of $F$. Let $s_w=\bigvee_{A\in\vec A_w}s_A$ assert that the switches occur in a pattern appearing in $\vec
A_w$. Finally, define $p_w=p_{[w]}\wedge s_w$.

We now prove that this works. In any forcing extension of $W$, we know by Lemma \ref{Lemma.LatticeLabels} that exactly one $p_{[u]}$ is true. And any
forcing extension exhibits some pattern $A$ of switches being true, and this $A$ is assigned to exactly one $w\in[u]$, so exactly one $s_w$ is true for
$w\in[u]$. Thus, in any forcing extension, exactly one $p_w=p_{[u]}\wedge s_w$ is true. We arranged that $p_{w_0}$ is true in $W$ by the assignment of
the pattern of switches holding in $W$ to the world $w_0$.

If $W^\P$ is a forcing extension satisfying $p_w$, then both $p_{[w]}$ and $s_w$ hold in $W^\P$. If $w\leq v$ in $F$, then we already know that
$p_{[v]}$ is forceable in $W^\P$, and $s_v$ is forceable from any extension of $W$ without affecting the buttons, so $p_v=p_{[v]}\wedge s_v$ is
forceable over $W^\P$, as desired. Conversely, if $p_v$ is forceable over $W^\P$, this implies that $p_{[v]}$ is forceable over $W^\P$, and so
$[w]\leq[v]$ in the quotient lattice. It follows that $w\leq v$ in $F$.\QED

We now state some definitions from modal logic. A {\df propositional world}, also called a {\df state}, is a map of the propositional variables to the
set $\{$true, false$\}$. This is simply a row in a truth table. A {\df Kripke model} $M$ is a set $U$ of propositional worlds, together with a relation
$R$ on $U$ called the {\df accessibility} relation. The Kripke semantics define when a modal assertion $\varphi$ is true at a world $w$ in a Kripke
model $M$, written $(M,w)\satisfies\varphi$. Namely, for atomic assertions, $(M,w)\satisfies q$ if $q$ is true in $w$; for Boolean connectives, one
uses the usual inductive treatment; for necessity, $(M,w)\satisfies\necessary\varphi$ if whenever $w\mathrel{R}v$, then $(M,v)\satisfies\varphi$; and
for possibility, $(M,w)\satisfies\possible\varphi$ if there is $v$ with $w\mathrel{R}v$ and $(M,v)\satisfies\varphi$. The underlying {\df frame} of the
model $M$ is the structure $\<U,R>$, ignoring the internal structure of the elements of $U$. The reader may easily check that every Kripke model whose
frame is a partial pre-order satisfies \S4, and every Kripke model on a directed partial pre-order satisfies \S{4.2}. A deeper fact is Lemma
\ref{Lemma.NotS4.2FailsOnALattice}, that the finite pre-lattice frames are complete for \S{4.2}. The next lemma is the heart of our argument, where we
prove that the behavior of any Kripke model on finite pre-lattice can be exactly simulated by forcing.

\SubLemma. If\/ $M$ is a Kripke model whose frame is a finite pre-lattice with a world $w_0$ and $W$ is a model of set theory with a sufficiently large
independent family of buttons and switches, then there is an assignment of the propositional variables $q_i$ to set theoretical assertions $\psi_i$,
such that for any modal assertion $\varphi$ we have \label{Lemma.InterpretingAnyKripkeModel}%
$$(M,w_0)\satisfies \varphi(q_0,\ldots,q_n)\quad\hbox{if and only if}\quad W\satisfies\varphi(\psi_0,\ldots,\psi_n).$$

\Proof. Each $\psi_i$ will be a certain Boolean combination of the buttons and switches. We have assumed that the frame $F$ of $M$ is a finite
pre-lattice. We may assume without loss of generality that $w_0$ is an initial world of $M$, by ignoring the worlds not accessible from $w_0$. Let
$p_w$ be the assertions assigned according to Lemma \ref{Lemma.PreLatticeLabels}. Since $w_0$ is an initial world of $F$, we may ensure that
$W\satisfies p_{w_0}$. Let $\psi_i=\bigvee\set{p_w\st (M,w)\satisfies q_i}$. We will prove the lemma by establishing the following stronger claim.
$$(M,w)\satisfies \varphi(q_0,\ldots,q_n)\quad\hbox{if and only if}\quad W\satisfies\necessary\left(\strut p_w\implies\varphi(\psi_0,\ldots,\psi_n)\right)$$
This is true for atomic $\varphi$, since $q_i$ is true at $w$ if and only if $p_w$ is one of the disjuncts of $\psi_i$, in which case
$p_w\implies\psi_i$ in any forcing extension of $W$, and conversely if $p_w\implies\psi_i$ is true in a forcing extension where $p_w$ is true, then
$\psi_i$ must be true there, in which case $q_i$ is true at $w$ in $M$. If the statement is true for $\varphi_0$ and $\varphi_1$, then it is also true
for $\varphi_0\wedge\varphi_1$. For negation, suppose that $(M,w)\satisfies\neg\varphi(q_0,\ldots,q_n)$. By induction,
$W\not\satisfies\necessary(\strut p_w\implies\varphi(\psi_0,\ldots,\psi_n))$, so there is a forcing extension $W^\P$ satisfying $p_w$ and
$\neg\varphi(\psi_0,\ldots,\psi_n)$. Since the truth values of $\psi_i$ necessarily depend only on the values of the various $p_u$, it follows that all
forcing extensions with $p_w$ will satisfy $\neg\varphi(\psi_0,\ldots,\psi_n)$. So we have proved $W\satisfies\necessary(\strut
p_w\implies\neg\varphi(\psi_0,\ldots,\psi_n))$, as desired, and reversing the steps establishes the converse. Finally,
$(M,w)\satisfies\possible\varphi(q_0,\ldots,q_n)$ if and only if $\exists u\geq w\,(M,u)\satisfies\varphi(q_0,\ldots,q_n)$, which occurs if and only if
$W\satisfies \necessary(\strut p_u\implies\varphi(\psi_0,\ldots,\psi_n))$. Since $W\satisfies\necessary(p_w\implies\possible p_u)$, this implies
$W\satisfies\necessary(\strut p_w\implies\possible\varphi(\psi_0,\ldots,\psi_n))$, as desired. Conversely, if $W\satisfies\necessary(\strut
p_w\implies\possible\varphi(\psi_0,\ldots,\psi_n))$, then $\varphi(\psi_0,\ldots,\psi_n)$ is forceable over any extension of $W$ with $p_w$. Since all
such extensions have $p_u$ for some $u\geq w$ and the $\psi_i$ depend only on the values of $p_v$, it must be that
$W\satisfies\necessary(p_u\implies\varphi(\psi_0,\ldots,\psi_n)$ for some $u\geq w$. By induction, this is equivalent to
$(M,u)\satisfies\varphi(q_0,\ldots,q_n)$ and consequently to $(M,w)\satisfies\possible\varphi(q_0,\ldots,q_n)$, as desired.\QED

The next step of our proof relies on a fact about \S{4.2}. A {\df tree} is a partial order $\leq$ on a set $F$ such that the predecessors of any node
are linearly ordered. A {\df pre-tree} is a partial pre-order $\leq$ on a set $F$ such that the quotient $F/{\equiv}$ is a tree; each node of this tree
is effectively replaced in $F$ with a cluster of equivalent nodes. A {\df baled tree} is a partial order $\leq$ on a set $F$ having a largest node
$b\in F$, such that $F\minus\singleton{b}$ is a tree (imagine baling or tying the top branches of a tree together, as in the figure below). A {\df
baled pre-tree} is the result of replacing each node in a baled tree with a cluster of equivalent nodes; equivalently, it is a partial pre-order whose
quotient by $\equiv$ is a baled tree. Note that every baled tree is a lattice, and every baled pre-tree is a pre-lattice. A partial pre-order is {\df
directed} if any two nodes have a common upper bound. A modal logic $\Theta$ has the {\df finite frame} property if whenever
$\Theta\not\proves\varphi$, then there is a finite $\Theta$-frame $F$ and a Kripke model having frame $F$ in which $\varphi$ fails.

\SubLemma. If a modal assertion $\varphi$ is not in \S{4.2}, then it fails in some Kripke model $M$ whose frame is a finite baled pre-tree, and hence a
finite pre-lattice.\label{Lemma.NotS4.2FailsOnALattice}

\Proof. It is easy to see that every Kripke model whose frame is a directed partial pre-order satisfies \S{4.2} and conversely that any frame that is
not a directed partial pre-order has a Kripke model violating \S{4.2} (see \cite[Theorem 3.38]{ChagrovZakharyaschev1997:ModalLogic}). This is what it
means to say that \S{4.2} is {\df defined} by the class of directed partial pre-orders. By \cite[Theorem 5.33]{ChagrovZakharyaschev1997:ModalLogic}, it
is known that \S{4.2} has the finite frame property, and so if $\S{4.2}\not\proves\varphi$, then $\varphi$ fails in a Kripke model $M_0$ whose frame
$F_0$ is a finite \S{4.2} frame, which is to say, a finite directed partial pre-order. We will construct a Kripke model $M$ that is bisimilar with
$M_0$, and which consequently has the same modal theory, but whose frame is a finite baled pre-tree, and consequently a finite pre-lattice. Our
construction is a minor modification of the standard technique of {\df tree unravelling} as described in the proof of \cite[Theorem
2.19]{ChagrovZakharyaschev1997:ModalLogic}. We know that $\varphi$ fails at some world $w_0$ in $M_0$, and we may assume that $w_0$ is in the smallest
cluster of $F_0$. By directedness, $F_0$ has a largest cluster $[b]$. The quotient $F_0/{\equiv}$ is a finite directed partial order. For each $[u]\in
F_0/{\equiv}$, let us say that $t$ is a path from $[w_0]$ to $[u]$ in $F_0/{\equiv}$ if it is a maximal linearly ordered subset of the interval
$[[w_0],[u]]$ in $F_0/{\equiv}$. Such paths form a tree when ordered by end-extension. Let $F$ be the {\df partial unravelling} of $F_0$, except for
the largest cluster.
$$
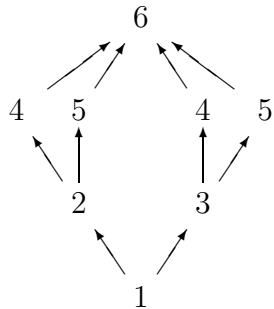
\hbox{\begin{diagram}[height=1.5em,width=1em]
            &            & \hbox{\tiny A partial order} &            &     &\\
            &            & 6 &            &     &\\
            &  \ruTo     &   & \luTo      &     &\\
 4          &            &   &            &  5  &\\
   \uTo     & \luTo(4,2) &   & \ruTo(4,2) & \uTo&\\
 2          &            &   &            &  3  &\\
            &   \luTo    &   & \ruTo      &     &\\
            &            & 1 &            &     &\\
\end{diagram}}\hskip3em
\hbox{\begin{diagram}[height=1.5em,width=1em,labelstyle=\scriptstyle,textflow]
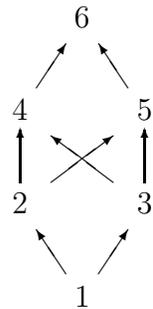
   &       &      &                 & \hbox{\tiny The resulting partial unravelling, a baled tree} &     &         &     &       \\
   &       &      &                 & 6                   &                    &         &     &         \\
   &       &      & \ruTo\ruTo(4,2) &                     & \luTo\luTo(4,2) &         &     &        \\
 4 &       & 5    &                 &                     &                    &    4    &     &  5       \\
   & \luTo & \uTo &                 &                     &                    &  \uTo   & \ruTo &        \\
   &       & 2    &                 &                     &                    &    3    &     &        \\
   &       &      & \luTo           &                     &      \ruTo         &         &     &       \\
   &       &      &                 &     1               &                    &         &     &       \\
\end{diagram}}$$
That is, $F$ consists of the maximal cluster $[b]$ of $F_0$, together with the set of all pairs $\<u,t>$, where $t$ is a path from $[w_0]$ to $[u]$ in
$F_0/{\equiv}$ and $[u]\neq[b]$. The order on $F$ is by end-extension of the paths $t$ and the $F_0$ order on $u$, with $[b]$ still maximal. The worlds
within any copy of a cluster are still equivalent and consequently still form a cluster, and so $F$ is a baled pre-tree. Let $M$ be the resulting
Kripke model on $F$, obtained by also copying the propositional values from every world $u\in F_0$ in $M_0$ to the copies $\<u,t>$ of it in $F$. It is
easy to see that $M$ is bisimilar with $M_0$, according to the correspondence that we have defined, because every world accesses in $M$ all and only
the copies of the worlds that it accesses in $M_0$. It follows that every world in $M_0$ satisfies exactly the same modal truths in $M_0$ that its
copies satisfy in $M$. Consequently, $\varphi$ fails at the copy of $w_0$ in $M$. Thus, $\varphi$ fails in a Kripke model whose frame is a finite baled
pre-tree, and all such frames are pre-lattices.\QED

\Proof Proof of Theorem \ref{MainTheorem}. Finally, we prove the theorem. By Theorem \ref{S4.2IsValid}, the \ZFC-provable principles of forcing
includes \S{4.2}. If $\varphi$ is not in \S{4.2}, then by Lemma \ref{Lemma.NotS4.2FailsOnALattice}, there is a Kripke model $M$ on a finite pre-lattice
in which $\varphi$ fails at an initial world. It is well known that if \ZFC\ is consistent, then so is $\ZFC+{\mathrm{V=L}}$, and so by Lemma
\ref{Lemma.ButtonsExist}, there is a model of set theory $\mathrm{L}$ having an infinite independent family of buttons and switches. By Lemma
\ref{Lemma.InterpretingAnyKripkeModel}, there is an assignment of the propositional variables of $\varphi$ to sentences $\psi_i$ such that
$\mathrm{L}\satisfies\neg\varphi(\psi_0,\ldots,\psi_n)$. Therefore, $\varphi$ is not a valid principle of forcing in $\mathrm{L}$, and hence not a
\ZFC-provable principle of forcing.\QED

\Section The Jankov-Fine formula \label{Section.TheJankovFineFormula}

While the previous section was written with a set theoretical reader in mind, let us now emphasize certain points for the modal logicians. The main
theorem can be restated in a way (as follows) that aligns it with many other completeness theorems in modal logic.

\Theorem. If\/ \ZFC\ is consistent, then $$\S{4.2}\proves\varphi(q_0,\ldots q_n)\quad\hbox{if and only if}\quad\forall
\psi_0,\ldots,\psi_n\,\,\ZFC\proves\varphi(\psi_0,\ldots,\psi_n),$$ where the $\psi_i$ range over the sentences in the language of set theory and
$\necessary$ and $\possible$ are understood in \ZFC\ with the forcing interpretation.\label{Theorem.MainTheoremRestatement}

Modal logicians will recognize that Lemmas \ref{Lemma.LatticeLabels} and \ref{Lemma.PreLatticeLabels} assert exactly that the assertions $p_w$ fulfill
the relevant Jankov-Fine formula, which we now define. For any graph $F=\<U,E>$, assign a propositional variable $p_w$ to each vertex $w\in U$ of the
graph and let $\delta(F)$ be the following formula, the {\df Jankov-Fine formula}. It asserts that, necessarily, exactly one $p_w$ is true, and if
$p_w$ is true, then $\possible p_v$ if and only if $w\mathrel{E}v$.
$$\begin{array}{clc}
\delta(F)= & \necessary\bigvee_{w\in U} p_w    &  \wedge\\
   & \necessary\bigwedge_{w\neq v}(p_w\implies\neg p_v)    &  \wedge\\
   & \necessary\bigwedge_{w\,E\,v}(p_w\implies\possible p_v)    &  \wedge\\
   & \necessary\bigwedge_{\neg w\,E\,v}(p_w\implies\neg\possible p_v)   &
\end{array}$$
Lemma \ref{Lemma.InterpretingAnyKripkeModel} has nothing essentially to do with set theory, but rather only with Kripke models and the Jankov-Fine
formula:

\SubLemma. Suppose that $M$ is a Kripke model whose frame $F$ is a finite partial pre-order and $w_0$ is a world of $M$. If\/ $N$ is any other Kripke
model satisfying \S4 (at some world $u_0$) and the Jankov-Fine formula $\delta(F)\wedge p_{w_0}$, then there is an assignment of the propositional
variables $q_i$ of $M$ to assertions $\psi_i$ in $N$ such that for any modal assertion $\varphi$,
\label{Lemma.IntrepretingKripkeModelsIfJankovFineIsSatisfied}
$$(M,w_0)\satisfies\varphi(q_0,\ldots,q_n)\quad\hbox{if and only if}\quad (N,u_0)\satisfies\varphi(\psi_0,\ldots,\psi_n)$$

\Proof. We argue just as in Lemma \ref{Lemma.InterpretingAnyKripkeModel}. Let $p_w$ be the assertions in $N$ satisfying the Jankov-Fine formula. Let
$\psi_i=\bigvee\set{p_w\st (N,w)\satisfies q_i}$. We establish the following stronger claim by induction on $\varphi$:
$$(M,w)\satisfies\varphi(q_0,\ldots,q_n)\quad\Leftrightarrow\quad (N,u_0)\satisfies\necessary(\strut p_w\implies \varphi(\psi_0,\ldots,\psi_n)).$$ The atomic case holds by the definition of $\psi_i$. Conjunction follows
because $\necessary$ distributes over $\wedge$. Negation follows via the properties of the Jankov-Fine formula, because every world in $N$ (accessible
from $u_0$) satisfies exactly one $p_w$, and any two such worlds agree on every $\psi_i$. Possibility follows using the Jankov-Fine formula again,
since $p_u$ is possible from a world with $p_w$ if and only if $w\leq u$ in $F$. Finally, the stronger claim implies the lemma, because $p_{w_0}$ is
true at $u_0$ in $N$.\QED

A modal theory $\Lambda$ is {\df closed under substitution} if $\varphi(\psi_0,\ldots,\psi_n)$ is in $\Lambda$ whenever $\varphi(q_0,\ldots,q_n)$ is,
for any modal assertions $\psi_i$. If every Kripke model with frame $F$ satisfies $\Lambda$, then $F$ is a {\df $\Lambda$-frame}.

\SubLemma. Suppose that $F$ is a finite partial pre-order. If a modal theory $\Lambda\fo\S4$ is closed under substitution and is consistent, for any
$w_0\in F$, with the Jankov-Fine formula $\delta(F)\wedge p_{w_0}$, then $F$ is a $\Lambda$-frame.\label{Lemma.MainModal}

\Proof. If $F$ is not a $\Lambda$-frame, then there is a Kripke model $M$ with frame $F$ and some $\varphi\in\Lambda$ such that
$(M,w_0)\satisfies\neg\varphi(q_0,\ldots,q_n)$ for some $w_0\in F$. If $\Lambda$ is consistent with $\delta(F)\wedge p_{w_0}$, then there is a Kripke
model satisfying $(N,u_0)\satisfies\Lambda\wedge\delta(F)\wedge p_{w_0}$. By Lemma \ref{Lemma.IntrepretingKripkeModelsIfJankovFineIsSatisfied}, there
is an assignment $q_i\mapsto\psi_i$, where $\psi_i$ is a Boolean combination of the $p_w$ in $N$, such that
$(N,u_0)\satisfies\neg\varphi(\psi_0,\ldots,\psi_n)$. This contradicts a substitution instance of $(N,u_0)\satisfies\Lambda$, since
$\varphi\in\Lambda$.\QED

Similarly, the concepts of button and switch are not set theoretic; they make sense in any Kripke model. Specifically, a {\df button} is a statement
that is necessarily possibly necessary, and a {\df switch} is a statement such that it and its negation are necessarily possible. A family of buttons
$\singleton{b_i}_{i\in I}$ and switches $\singleton{s_j}_{j\in J}$ is {\df independent} in $M$ at world $u$ if none of the buttons is necessary at $u$
and necessarily, any button can be turned on and any switch can be turned either on or off without affecting the other buttons and switches. This can
be expressed precisely in modal logic as follows. For any $A\of I$ and $B\of J$, let $\Theta_{A,B}=(\bigwedge_{i\in A}\necessary
b_i)\wedge(\bigwedge_{i\notin A}\neg\necessary b_i)\wedge(\bigwedge_{j\in B}s_j)\wedge(\bigwedge_{j\notin B}\neg s_j)$ assert that the pattern of
buttons and switches is specified by $A$ and $B$. A family $\singleton{b_i}_{i\in I}\union \singleton{s_j}_{j\in J}$ of buttons and switches is
independent if:
$$\bigl(\bigwedge_{i\in I}\neg\necessary b_i\bigr)\wedge
  \bigwedge_{\vtop{\tiny\hbox{$A\of I$}\hbox{$B\of J$}}}
     \necessary\bigl(\Theta_{A,B}\implies\bigwedge_{\vtop{\tiny\hbox{$A\of A'$}\hbox{$B'\of J$}}}\possible\Theta_{A',B'}\bigr).$$
Thus, the buttons are off initially, and necessarily, from any possible pattern of buttons and switches, any larger pattern of buttons and any pattern
of switches is possible. The main technique in our proofs of Lemmas \ref{Lemma.LatticeLabels} and \ref{Lemma.PreLatticeLabels} appears to be very
reminiscent of Smory\'{n}ski's \cite{Smorynski1993:ApplicationsOfKripkeModels} proof of de Jongh's theorem \cite{deJongh1970:MaximalityIPCwrtHA} on
Heyting's Arithmetic.

\SubLemma. If\/ $F$ is a finite pre-lattice, $w_0\in F$ and $\Lambda\fo\S4$ is a modal theory consistent with a sufficiently large independent family
of buttons and switches, then $\Lambda$ is consistent with the Jankov-Fine formula $\delta(F)\wedge
p_{w_0}$.\label{Lemma.IfConsistentWithButtonsThenJankovFine}

\Proof. Suppose that $(M,u_0)\satisfies\Lambda$ has a sufficiently large independent family of buttons and switches. The proofs of Lemmas
\ref{Lemma.LatticeLabels} and \ref{Lemma.PreLatticeLabels} work in $M$. Specifically, in those arguments we assigned to each node $w$ in $F$ an
assertion $p_w$, a Boolean combination of buttons and switches, so that at any world $u$ accessible from $u_0$ in $M$, exactly one of the $p_w$ is
true, and if $p_w$ is true then $\possible p_v$ holds if and only if $u\leq v$ in $F$. Thus, $(M,u_0)\satisfies\delta(F)$. By assigning the pattern of
switches that happens to hold in $M$ at $u_0$ to the node $w_0$, we also arranged $(M,u_0)\satisfies p_{w_0}$. So $\Lambda$ is consistent with the
Jankov-Fine formula $\delta(F)\wedge p_{w_0}$.\QED

The lemmas combine to prove Theorems \ref{MainTheorem} and \ref{Theorem.MainTheoremRestatement} as follows. Let $\Lambda$ be the \ZFC-provable
principles of forcing. It is easy to see that $\Lambda$ is closed under substitution, modus ponens and necessitation. By Theorem \ref{S4.2IsValid}, we
know $\S{4.2}\of\Lambda$. By Lemma \ref{Lemma.ButtonsExist}, if \ZFC\ is consistent, then there are models of set theory having infinite independent
families of buttons and switches. It follows that $\Lambda$ is consistent with arbitrarily large finite independent families of buttons and switches.
By Lemma \ref{Lemma.IfConsistentWithButtonsThenJankovFine}, therefore, $\Lambda$ is consistent with the Jankov-Fine formula $\delta(F)\wedge p_{w_0}$
for any finite pre-lattice $F$. By Lemma \ref{Lemma.MainModal}, therefore, all such $F$ are $\Lambda$-frames. By Lemma
\ref{Lemma.NotS4.2FailsOnALattice}, any statement not in \S{4.2} fails in a Kripke model having such a frame and consequently is not in $\Lambda$. So
$\Lambda\of\S{4.2}$ and, consequently, $\Lambda=\S{4.2}$.

Let us now push these techniques a bit harder, in order to arrive at a new class of frames complete for \S{4.2} and some useful characterizations of
\S4, \S{4.2} and \S5. A partial pre-order $(B,\leq)$ is a {\df pre-Boolean algebra} if the quotient partial order $B/{\equiv}$ is a Boolean algebra.

\Lemma. For any natural numbers $n$ and $m$, there is a Kripke model $N$ whose frame is a finite pre-Boolean algebra, such that at any initial world in
$N$, there is an independent family of $n$ buttons and $m$ switches.\label{Lemma.PreBooleanAlgebrasHaveButtonsSwitches}

\Proof. Using power sets, let $F=\mathrm{P}(n)\cross\mathrm{P}(m)$, so that the nodes of $F$ consist of pairs $(A,B)$, where $A\of
n=\set{0,\ldots,n-1}$ and $B\of m=\set{0,\ldots,m-1}$. The order is determined by the first coordinate only, so that $(A,B)\leq (A',B')$ if and only if
$A\of A'$. This is clearly a partial pre-order. The corresponding equivalence relation is $(A,B)\equiv (A',B')$ if and only if $A=A'$, and so the
quotient is isomorphic to the power set $\mathrm{P}(n)$, which is a finite Boolean algebra. So $F$ is a finite pre-Boolean algebra. Let $N$ be the
Kripke model on $F$ in which $b_i$ true at $(A,B)$ when $i\in B$ and $s_j$ true at $(A,B)$ when $j\in B$. Clearly, every $b_i$ is a button in $N$ and
every $s_j$ is a switch in $N$, and they are independent at any initial world of $N$, because whatever the pattern $(A,B)$ of buttons and switches in
any world of $N$, any larger pattern of buttons $A'\fo A$ and any pattern of switches $B'\of m$ is possible.\QED

The Kripke model $N$ produced in Lemma \ref{Lemma.PreBooleanAlgebrasHaveButtonsSwitches} has the smallest frame supporting an independent family of $n$
buttons and $m$ switches, because for independence one one needs worlds realizing every pattern $(A,B)$ of buttons and switches. A class $\F$ of frames
is {\df complete} for a modal theory $\Lambda$ if every $F\in{\F}$ is a $\Lambda$-frame and any $\varphi$ true in all Kripke models having frames in
$\F$ is in $\Lambda$.

\Lemma. A class $\F$ of frames is complete for \S{4.2} if and only if every $F\in\F$ is an \S{4.2} frame and there are Kripke models, with frames in
$\F$, for arbitrarily large finite independent families of buttons and switches.\label{Lemma.CompleteForS4.2IffButtonsAndSwitches}

\Proof. The forward implication is immediate, because \S{4.2} is consistent by Lemma \ref{Lemma.PreBooleanAlgebrasHaveButtonsSwitches} with the
existence of large independent families of buttons and switches. Conversely, suppose the latter property. If $\varphi_0$ is not in \S{4.2}, then by
Lemma \ref{Lemma.NotS4.2FailsOnALattice} there is some Kripke model $M$ whose frame $F$ is a finite pre-lattice such that
$(M,w_0)\satisfies\neg\varphi_0(q_0,\ldots,q_n)$. Let $(N,u_0)$ be a Kripke model with frame in $\F$ having an independent family of $n$ buttons and
$m$ switches, where $n$ is the number of clusters in $F$ and the size of any cluster is at most $2^m$. By the proof of Lemma
\ref{Lemma.IfConsistentWithButtonsThenJankovFine}, there are assertions $p_w$ for $w\in F$ such that $(N,u_0)\satisfies\delta(F)\wedge p_{w_0}$. By
Lemma \ref{Lemma.IntrepretingKripkeModelsIfJankovFineIsSatisfied}, there is an assignment of the propositional variables $q_i$ of $M$ to assertions
$\psi_i$ in $N$ such that $(M,w_0)\satisfies\varphi(q_0,\ldots,q_n)$ if and only if $(N,u_0)\satisfies\varphi(\psi_0,\ldots,\psi_n)$. By the assumption
on $\varphi_0$, this means that $(N,u_0)\satisfies\neg\varphi_0(\psi_0,\ldots,\psi_n)$. Thus, a substitution instance of $\varphi_0$ fails at a world
in $N$, a Kripke model whose frame is in $\F$. So this class of frames is complete for \S{4.2}.\QED

\Theorem. The following sets of frames are complete for \S{4.2}.
\begin{enumerate}
 \item Finite directed partial pre-orders.
 \item Finite pre-lattices.
 \item Finite baled pre-trees.
 \item Finite pre-Boolean algebras.
\end{enumerate}

\Proof. All of these frames are directed partial pre-orders, and so they are all \S{4.2} frames. Lemma \ref{Lemma.NotS4.2FailsOnALattice} shows that
any statement not in \S{4.2} fails in a Kripke model whose frame is a finite baled pre-tree, and hence a finite pre-lattice and a finite pre-order, so
these classes are complete for \S{4.2}. The new part of this theorem is (4). By Lemma \ref{Lemma.PreBooleanAlgebrasHaveButtonsSwitches} the class of
finite pre-Boolean algebras have Kripke models for arbitrarily large finite independent families of buttons and switches. So by Lemma
\ref{Lemma.CompleteForS4.2IffButtonsAndSwitches}, this class is also complete for \S{4.2}.\QED

Let us summarize what we have proved about \S{4.2}:

\Theorem. Suppose that the modal theory $\Lambda$ contains \S4 and is closed under substitution. For any class $\F$ of finite frames complete for
\S{4.2}, the following are equivalent:
\begin{enumerate}
 \item $\Lambda$ is consistent with arbitrarily large finite independent families of buttons and switches.
 \item $\Lambda$ is consistent with the Jankov-Fine formula $\delta(F)\wedge p_w$ for any frame $F\in\F$ and world $w\in F$.
 \item Every frame in $\F$ is a $\Lambda$-frame.
 \item $\Lambda\of\S{4.2}$.
\end{enumerate}\label{Theorem.EquivalentsForBeingInsideS4.2}

\Proof. Lemma \ref{Lemma.MainModal} shows that (2) implies (3), as the elements of $\F$ must be finite directed partial pre-orders. For (3) implies
(4), observe that if $\varphi\notin\S{4.2}$, then it must fail in a Kripke model whose frame is in $\F$, contrary to (3). For (4) implies (2), note
that \S{4.2} and $\delta(F)\wedge p_{w_0}$ are true together in the Kripke model having frame $F$, with $p_w$ true exactly in world $w$. For (4)
implies (1), if $\Lambda\of\S{4.2}$, then $\Lambda$ is true in the Kripke models constructed in Lemma
\ref{Lemma.PreBooleanAlgebrasHaveButtonsSwitches}, which have large independent families of buttons and switches. Finally, Lemma
\ref{Lemma.IfConsistentWithButtonsThenJankovFine} shows that (1) implies (2) in the special case where $\F$ is the class of all finite pre-lattices,
and hence (1) implies (4) for any $\F$, since they do not depend on $\F$, completing the proof.\QED

Our later analysis will benefit from similar characterizations of \S4 and \S5. For \S5, we use the fact that the class of finite complete reflexive
graphs are complete for \S5. This result will be applied in Theorem \ref{Theorem.ValiditiesAreWithinS5}.

\Theorem. Suppose that a modal theory $\Lambda$ contains \S4 and is closed under substitution. Then the following are equivalent:
\begin{enumerate}
 \item $\Lambda$ is consistent with arbitrarily large finite independent families of switches.
 \item $\Lambda$ is consistent with the Jankov-Fine formula $\delta(F)\wedge p_w$ for any finite complete graph $F$ and world $w\in F$.
 \item Every finite complete reflexive graph is a $\Lambda$-frame.
 \item $\Lambda\of\S5$.
\end{enumerate}\label{Theorem.EquivalentsForBeingInsideS5}

\Proof. For (1) implies (2), the point is that when $F$ has only one cluster, the argument of Lemma \ref{Lemma.PreLatticeLabels} does not require any
buttons. Suppose that $N$ is a Kripke model having an independent family of switches $\singleton{s_j}_{j\in J}$. For any $A\of J$, define $s_A$ as in
Lemma \ref{Lemma.PreLatticeLabels} to assert that the pattern of switches is $A$. Partition the collection of $A\of J$ among the worlds $u\in F$ by
assigning a nonempty set $\vec A_u$ of sets to each world $u$ in $F$. Let $s_u=\bigwedge_{A\in\vec A_u}s_A$ as in Lemma \ref{Lemma.PreLatticeLabels}.
By assigning whatever pattern of switches holds at $u_0$ to the world $w$, we can arrange that $(N,u_0)\satisfies s_w$. Since every world must have
some unique pattern of switches, it follows that $(M,u_0)$ satisfies that necessarily, exactly one $s_u$ is true. And since the switches are
independent, we also know that $(M,u_0)\satisfies\necessary s_u$ for any $u\in F$. Since all worlds in $F$ are accessible from each other, this implies
$(N,u_0)\satisfies\delta(F)\wedge p_w$, as desired.

Lemma \ref{Lemma.MainModal} shows that (2) implies (3). And (3) implies (4) because any statement not in \S5 is known to fail in a Kripke model whose
frame is a finite complete reflexive graph (see \cite[Proposition 3.32, Corollary 5.29]{ChagrovZakharyaschev1997:ModalLogic}). Finally, if $\Lambda\of
\S5$, then $\Lambda$ holds in any Kripke model whose frame is a complete graph. It is easy to arrange independent families of switches in such Kripke
models, just by ensuring that every possible pattern of switches is exhibited in some world.\QED

In the \S4 context, one can generalize Theorems \ref{Theorem.EquivalentsForBeingInsideS4.2} and \ref{Theorem.EquivalentsForBeingInsideS5} to the
following, which we expect will be relevant for Conjecture \ref{Conjecture.CCCvaliditiesAreS4}. Recall that a modal logic $\Theta$ has the finite frame
property if whenever $\Theta\not\proves\varphi$, then there is a finite $\Theta$-frame $F$ and a Kripke model $M$ with frame $F$ in which $\varphi$
fails. Such theories as \S4, \S{4.2}, \S{4.3} and \S5 are known by \cite[5.29, 5.32, 5.33]{ChagrovZakharyaschev1997:ModalLogic} to have the finite
frame property.

\Theorem. Suppose that $\Theta$ is a modal logic containing \S4 and has the finite frame property, such as \S4, \S{4.2}, \S{4.3} or \S5. If\/ $\Lambda$
contains \S4 and is closed under substitution, then for any class $\F$ of frames complete for $\Theta$, the following are equivalent:
\begin{enumerate}
 \item $\Lambda$ is consistent with the Jankov-Fine formula $\delta(F)\wedge p_{w_0}$ for any frame $F\in\F$ and world $w_0\in F$.
 \item Every frame in $\F$ is a $\Lambda$-frame.
 \item $\Lambda\of\Theta$.
\end{enumerate}\label{Theorem.CharacterizationOfS4+S4.2+S5}

\Proof. If the Jankov-Fine formula $\delta(F)\wedge p_{w_0}$ holds in a Kripke model $(N,u_0)\satisfies\Lambda$, then by Lemma
\ref{Lemma.IntrepretingKripkeModelsIfJankovFineIsSatisfied}, any Kripke model $M$ with frame $F$ has an assignment $q_i\mapsto\psi_i$ such that
$(M,w_0)\satisfies\varphi(q_0,\ldots,q_n)$ if and only if $(N,u_0)\satisfies\varphi(\psi_0,\ldots,\psi_n)$. In particular, if $\varphi$ fails in such
an $M$, then the corresponding substitution instance of it will fail in $N\satisfies\Lambda$. So every such $F$ is a $\Lambda$-frame, and we have
proved (1) implies (2). For (2) implies (3), observe that if $\varphi$ is not in $\Theta$, then it fails in a Kripke model $M$ having a finite frame
$F\in\F$. Since $F$ is a $\Lambda$-frame, we know $M\satisfies\Lambda$ and consequently $\varphi$ is not in $\Lambda$. So $\Lambda\of\Theta$. For (3)
implies (1), observe that the Jankov-Fine formula $\delta(F)\wedge p_{w_0}$ is easily satisfied at a Kripke model $M$ having frame $F$, where $p_w$ is
true exactly at $w$. Since $F$ is a $\Theta$-frame, we know $M\satisfies\Theta$ and consequently $M\satisfies\Lambda$ by (3), so (1) holds, and the
proof is complete.\QED

We do not have a button-and-switch characterization of \S4 in Theorem \ref{Theorem.CharacterizationOfS4+S4.2+S5}, because the frame of an \S4 model
need not be directed, and it is not true that every possibly necessary statement is necessarily possibly necessary. Under \S4, one can have assertions
$\varphi$ such that, simultaneously, $\varphi$ is possibly necessary, $\neg\varphi$ is possibly necessary, and $\varphi$ is possibly a switch. Thus,
for the \S4 context, we emphasize the official definition of button as a statement that is {\it necessarily} possibly necessary. Such examples show
that unlike \S{4.2}, under \S4 it is no longer true that every statement is a button, the negation of a button or a switch.

\Section Forcing over a fixed model of set theory \label{Section.ForcingOverFixedModel}

While our main theorem establishes that the \ZFC-provable principles of forcing are exactly those in \S{4.2}, it is not true that every model of set
theory observes only these validities. For any $W\satisfies\ZFC$, recall that a modal assertion $\varphi$ is a {\df valid principle of forcing in $W$}
if for all sentences $\psi_i$ in the language of set theory we have $W\satisfies\varphi(\psi_0,\ldots,\psi_n)$. For meta-mathematical reasons connected
with Tarski's theorem on the non-definability of truth, there is initially little reason to expect that the collection of such $\varphi$ should be
definable in $W$; rather, the assertion that $\varphi$ is valid in $W$ is expressed as a scheme, asserting all substitution instances
$\varphi(\psi_0,\ldots,\psi_n)$ in $W$. So this is formally a second-order notion. Because a statement is provable in \ZFC\ exactly when it holds in
all models of \ZFC, our main theorem establishes that the modal assertions $\varphi$ that are valid in all models of set theory are exactly those in
\S{4.2}. Our proof, however, established the following stronger result.

\Theorem. If\/ $W$ is a model of set theory with arbitrarily large finite independent families of buttons and switches, then the valid principles of
forcing in $W$ are exactly \S{4.2}.\label{Theorem.IfButtonsThenS4.2}

This result is stronger because it shows that the minimal set of forcing validities is realized in a single model of set theory (such as any model of
$\mathrm{V=L}$), rather than arising as the intersection of the validities of several models. Nevertheless, there are models of set theory whose valid
principles of forcing go beyond \S{4.2}. For example, the Maximality Principle \MP\ of \cite{Hamkins2003:MaximalityPrinciple} asserts all instances of
the scheme $\possible\necessary\psi\implies\psi$, so that any set theoretic statement $\psi$ that holds in some forcing extension and all further
extensions is already true. In other words, \MP\ asserts that \S5 is valid for forcing. Because it is established in
\cite{Hamkins2003:MaximalityPrinciple} that if \ZFC\ is consistent, then so is $\ZFC+\MP$, we conclude:

\Theorem. If\/ \ZFC\ is consistent, then it is consistent with \ZFC\ that all \S5 assertions are valid principles of forcing.

The forcing validities of a model, however, never go beyond \S5.

\Theorem. The valid principles of forcing in any model of set theory are included within \S5.\label{Theorem.ValiditiesAreWithinS5}

\Proof. Let $\Lambda$ be the set of forcing validities in a model $W$ of set theory. By Theorem \ref{S4.2IsValid}, this includes all of \S{4.2}. Also,
$\Lambda$ is easily seen to be closed under substitution. Observe next that any model of set theory has an infinite independent family of switches,
such as $s_n=$``the \GCH\ holds at $\aleph_n$.'' These and their negations are forceable in any finite pattern over any model of set theory by well
known forcing arguments. It follows that $\Lambda$ is consistent with the modal assertions that there are large independent families of switches. By
Theorem \ref{Theorem.EquivalentsForBeingInsideS5}, consequently, $\Lambda\of\S5$.\QED

\Corollary. If\/ \ZFC\ is consistent, then there is a model of set theory whose valid principles of forcing are exactly \S5.

Our results establish that both \S{4.2} and \S5 are realized as the exact set of forcing validities of a model of set theory (realized, respectively,
in models of $\mathrm{V=L}$ or of \MP).

\Question. Which modal theories arise as the valid principles of forcing in a model of \ZFC? \label{Question.WhichModalTheoriesArise?}

For example, can there be a model of set theory whose valid principles of forcing are exactly \theoryf{Dm.2} or exactly \S{4.3}? We have seen that if a
model of set theory has sufficiently many independent buttons and switches, then the valid principles of forcing will be only \S{4.2}. At the other
extreme, if there are no buttons, then \MP\ holds and so the valid principles of forcing in the model will be \S5. Is it possible to have a model with
a finite bound on the size of an independent family of buttons?

\Question. Is there a model of \ZFC\ with one unpushed button but not two independent buttons?

What are the valid principles of forcing in such models? It is clear that if there are only finitely many buttons in a model (meaning that the
independent families of buttons have some bounded finite size), then we could use Lemmas \ref{Lemma.PreLatticeLabels} and
\ref{Lemma.InterpretingAnyKripkeModel} to simulate any Kripke model on a pre-lattice frame with a correspondingly bounded number of clusters. The
various classes of such frames generate modal theories strictly between \S{4.2} and \S5. Are these realizable as the forcing validities of models of
set theory? Are these the only modal theories that arise? The next theorem is a start on these questions. Denote the collection of modal assertions
$\varphi$ that are valid in $W$ by $\Force^W$. Two buttons $b_0$ and $b_1$ are {\df semi-independent} if both are unpushed and one can push $b_0$
without pushing $b_1$.

\Theorem. Suppose that $W$ is a model of \ZFC\ set theory.
\begin{enumerate}
 \item $\S{4.2}\of\Force^W\of\S5$.
 \item If $W$ has no unpushed buttons, then $\Force^W=\S5$.
 \item If $W$ has at least one unpushed button, then $\Force^W\ofnoteq\S5$.
 \item If $W$ has two semi-independent buttons, then $\axiomf{W5}$ is not valid in $W$, and so $\S{4W5}\not\of\Force^W$.
 \item If $W$ has two independent buttons, then $\axiomf{.3}$ is not valid in $W$, and so $\S{4.3}\not\of\Force^W$.
 \item If $W$ has an independent family of one button and one switch, then $\axiomf{Dm}$ is not valid in $W$, and so $\theoryf{Dm}\not\of\Force^W$.
\end{enumerate}\label{Theorem.ButtonAndSwitchCharacterizationsOfForce^W}

\Proof. Statement (1) is the content of Theorems \ref{S4.2IsValid} and \ref{Theorem.ValiditiesAreWithinS5}. Over \S4, the additional \S5 axiom is
equivalent to $\possible\necessary\varphi\implies\necessary\varphi$, which exactly asserts that every button is pushed. So (2) and (3) hold. For (4),
suppose that $b_0$, $b_1$ are semi-independent buttons in $W$, and let $\varphi=(\neg\necessary b_0\wedge\neg\necessary b_1)\vee\necessary(b_0\wedge
b_1)$, which asserts that either none or both buttons are pushed. Thus, $\varphi$ is both true and possibly necessary in $W$, since the buttons are
initially unpushed and we could push both, but $\varphi$ is not necessary, since we could push just $b_0$; this violates \axiomf{W5}. For (5), suppose
that $W$ has two independent buttons $b_1$ and $b_2$. We argue as in Observation \ref{Observation.NotValidBeyondS4.2}. Let $\varphi=\necessary
b_1\wedge\neg\necessary b_2$ and $\psi=\necessary b_2\wedge\neg\necessary b_1$. Since both buttons are unpushed in $W$ and either may be pushed, we
conclude that $W\satisfies\possible\varphi\wedge\possible\psi$. But in any forcing extension of $W$, if $\varphi$ is true, then $\psi$ is impossible
and vice versa. So this instance of the conclusion of $.3$ fails in $W$. Thus, .3 is not valid in $W$, and so $\Force^W$ does not include \S{4.3}. For
(6), suppose that $W$ has an independent family of one button $b$ and one switch $s$. We follow the argument of Observation
\ref{Observation.NotValidBeyondS4.2}, case \axiomf{Dm}. We may assume that both $\necessary b$ and $s$ are false in $W$. Let $\varphi=\necessary b\vee
s$. This is possibly necessary in $W$, since one could push the button $b$, but not true in $W$, so $\possible\necessary\varphi\implies\varphi$ is
false in $W$. In any forcing extension of $W$, if $\necessary(\varphi\implies\necessary\varphi)$, then it must be that the button $b$ has been pushed
there, since otherwise one could have $\varphi$ true and then false again by flipping the switch $s$. So we have argued that
$(\necessary(\varphi\implies\necessary\varphi))\implies\varphi$ holds in every forcing extension of $W$. Thus, the hypothesis of this instance of
\axiomf{Dm} holds in $W$, while the conclusion fails, so \theoryf{Dm} is not valid in $W$. This establishes (6).\QED

Semi-independent buttons are part of the following more general arrangement. A list of assertions $\varphi_1,\ldots,\varphi_n$ is a {\df volume
control} if each is a button, pushing any of them necessarily also pushes all the previous, and any of them can be pushed without pushing the next.
More precisely:
$$\necessary\bigwedge_{i<n}(\necessary\varphi_{i+1}\implies\necessary\varphi_i)\wedge(\possible\necessary\varphi_{i+1})\wedge
\bigl(\neg\necessary\varphi_{i+1}\implies\possible(\neg\varphi_{i+1}\wedge\necessary\varphi_i)\bigr).$$ The idea is that one can turn up the volume to
level $j$ by forcing $\necessary\varphi_j$, but there is no turning it down again. The volume control has volume zero if $\neg\necessary\varphi_1$ (and
so a volume control of length $n$ has $n+1$ many volume settings). These volume controls exhibit dependence of buttons, in a linear fashion, rather
than independence. If buttons $b_0$ and $b_1$ are semi-independent, then $\necessary b_0$, $\necessary(b_0\wedge b_1)$ is a volume control of length 2;
conversely, every volume control of length $2$ (and volume zero) consists of two semi-independent buttons. Similar ideas with more buttons produce
arbitrarily long volume controls. A Kripke model with a linear frame of $n+1$ clusters admits volume controls of length $n$ but no independent buttons.
A family of volume controls, buttons and switches is {\df independent} in a model if all the volume controls have zero volume in that model, all the
buttons are unpushed in that model, and in any forcing extension, one can operate any of the volume controls, buttons and switches without affecting
any of the others.

\Theorem. If\/ $W$ is a model of set theory exhibiting arbitrarily long volume controls independent from arbitrarily large families of independent
switches, then $\Force^W\of\S{4.3}$.

\Proof. We know by \cite[3.31, 3.32, 3.40, 5.33]{ChagrovZakharyaschev1997:ModalLogic} that the finite linear pre-orders are a complete class of \S{4.3}
frames. Suppose that $F$ is such a finite linear pre-order and $w_0$ is a node in the minimal cluster of $F$. Let $v_1,\ldots,v_n$ be a volume control
in $W$, where $n$ is the number of clusters of $F$, and suppose that this volume control forms an independent family with the switches $s_0,\ldots,
s_m$, where the size of any cluster is at most $2^m$. As in Lemmas \ref{Lemma.LatticeLabels} and \ref{Lemma.PreLatticeLabels}, we will assign to each
node $w\in F$ an assertion $p_w$, so that $W$ satisfies the Jankov-Fine formula $\delta(F)\wedge p_{w_0}$. Specifically, within each cluster $[u]$,
assign to each node $w\in[u]$ a nonempty set $\vec A_w$ of subsets $A\of n$ in such a way that these partition all subsets $A\of n$, and define
$s_w=\bigvee_{A\in\vec A_w}s_A$, as in Lemma \ref{Lemma.PreLatticeLabels}. In the least cluster, we assign the pattern $A$ of switches holding in $W$
to the node $w_0$. Now, for any node $w\in F$, if $w$ is in the $i^\th$ cluster, then we define $p_w=v_i\wedge s_w$. That is, the volume control
indicates the intended cluster and the switch indicates the intended node within that cluster. Since every forcing extension of $W$ exhibits some
volume setting and some pattern of switches, it is clear that it will satisfy exactly one of the $p_w$, and in any extension of $W$ where $p_w$ holds,
then $\possible p_v$ holds if and only if $w\leq v$. Thus, the Jankov-Fine formula $\delta(F)\wedge p_{w_0}$ is satisfied. It now follows by Theorem
\ref{Theorem.CharacterizationOfS4+S4.2+S5} that $F$ is a $\Force^W$-frame. So every \S{4.3}-frame is a $\Force^W$-frame, and so
$\Force^W\of\S{4.3}$.\QED

We close this section with a curious question. A modal theory is {\df normal} if it is closed under modus ponens, substitution and necessitation. All
the named modal theories that we have considered, such as \S4, \S{4.2}, \S{4.3}, \theoryf{Dm.2}, \S5 and so on, are normal. Meanwhile, the valid
principles of forcing in any model of \ZFC\ is easily seen to be closed under modus ponens and substitution. But is it closed under necessitation?

\Question. If\/ $\varphi$ is a valid principle of forcing, does it remain valid in all forcing extensions?

In other words, is $\necessary\varphi$ also valid? Equivalently, if $W\satisfies\ZFC$, then is $\Force^W$ normal? Of course, if $\Force^W$ is \S{4.2}
or \S5, then the answer is yes, so this question is related to Questions \ref{Question.WhichModalTheoriesArise?} and
\ref{Question.CanForceBeBetweenS4.2AndS5?}.

\Section The modal logic of forcing with parameters \label{Section.WithParameters}

We know from \cite{Hamkins2003:MaximalityPrinciple} that parameters play a subtle role in the strength of the Maximality Principles. While \MP\ is
equiconsistent with \ZFC, allowing real parameters in the scheme results in the principle $\MP(\R)$, which has some large cardinal strength; allowing
uncountable parameters leads to inconsistency. Allowing real parameters from all forcing extensions leads to a principle $\necessary\MP(\R)$ with a
large cardinal strength of at least infinitely many Woodin cardinals (but less than $\AD_\R+\Theta$ is regular). So let us analyze the role played by
parameters in the valid principles of forcing. Specifically, define that $\varphi(q_0,\ldots,q_n)$ is a {\df valid principle of forcing in a model $W$
with parameters in $X$} if for any set theoretical formulas $\psi_i(\vec x)$ we have $$W\satisfies \forall \vec x\in X\ \varphi(\psi_0(\vec
x),\ldots,\psi_n(\vec x)).$$ We denote the collection of such $\varphi$ by $\Force^W(X)$. The next theorem provides another answer to Question
\ref{MainQuestion}.

\Theorem. In any model of set theory, the modal assertions $\varphi$ that are valid principles of forcing with all parameters are exactly those in
\S{4.2}. Succinctly, $\Force^W(W)=\S{4.2}$.\label{Theorem.ValidWithParametersIsS4.2}

\Proof. Certainly any \S{4.2} assertion is a valid principle of forcing with any parameters, by the argument of Theorem \ref{S4.2IsValid}, which did
not depend on whether there were parameters or not. For the converse direction, we argue as in the Main Theorem, using the following fact, that with
sufficient parameters one can always construct independent families of buttons and switches.

\SubLemma. If\/ $W$ is any model of \ZFC, then with parameters $\omega_n^W$ for $n\in\omega$, there is an infinite independent family of buttons and
switches.

\Proof. Let $b_n$ be the assertion ``$\omega_{n+1}^W$ is not a cardinal'' and let $s_m$ be ``the \GCH\ holds at $\aleph_{\omega+m}$,'' referring to
$\aleph_{\omega+m}$ {\it de dicto}, rather than with a parameter. In any forcing extension of $W$, any of the buttons can be forced without affecting
the truth of the other buttons, by collapsing $\omega_{n+1}^W$ to its predecessor (see \cite[15.21]{Jech:SetTheory3rdEdition}). After this, the
switches $s_m$ and their negations can be forced in any desired pattern without adding bounded sets below $\aleph_\omega$ and, consequently, without
affecting the buttons $b_n$ (see \cite[15.18 \& related]{Jech:SetTheory3rdEdition}).\QED

The result now follows from Theorem \ref{Theorem.IfButtonsThenS4.2} and the observation that the existence of parameters can simply be carried through
that argument. Specifically, if $\varphi$ is not in \S{4.2}, then it fails at some world $w_0$ in a Kripke model $M$ whose frame is a finite
pre-lattice. As in the Main Theorem, we use the buttons and switches to define $p_w$ for each world $w$ in $M$ verifying the Jankov-Fine formula as in
Lemma \ref{Lemma.PreLatticeLabels}. If we define $\psi_i$ as in Lemma \ref{Lemma.InterpretingAnyKripkeModel}, then we observe as before that
$(M,w_0)\satisfies\varphi(q_0,\ldots,q_n)$ if and only if $W\satisfies\varphi(\psi_0,\ldots,\psi_n)$. Since $\varphi$ fails at $w_0$ in $M$ and the
$\psi_i$ are Boolean combinations of the buttons and switches, this produces a failing substitution instance of $\varphi$ in $W$ using the same
parameters.\QED

\Theorem. For any $W\satisfies\ZFC$ and any set $X$ of parameters in $W$, $$\S{4.2}\of\Force^W(X)\of\Force^W\of\S5.$$ If $X\of Y$, then
$$\S{4.2}\of\Force^W(Y)\of\Force^W(X)\of\Force^W\of\S5.$$

\Proof. Certainly any \S{4.2} assertion is valid, even with parameters, so $\S{4.2}\of\Force^W(X)$. The proof of Theorem
\ref{Theorem.ValiditiesAreWithinS5} is not affected by the presence of parameters, so $\Force^W(X)\of\S5$. If $X\of Y$, then it is at least as hard for
a modal assertion $\varphi$ to be valid for all substitution instances using parameters in $Y$ as for parameters in $X$, so
$\Force^W(Y)\of\Force^W(X)$. Finally, $\Force^W=\Force^W(\emptyset)$.\QED

\Question. Can the set of forcing validities $\Force^W(X)$ be strictly between \S{4.2} and \S5? When is it equal to \S{4.2} or to
\S5?\label{Question.CanForceBeBetweenS4.2AndS5?}

For example, in the proof of Theorem \ref{Theorem.ValidWithParametersIsS4.2}, we only used parameters $\omega_n^W$, so we may conclude:

\Corollary. For any model $W$ of set theory, $\Force^W(\aleph_\omega^W)=\S{4.2}$. More specifically, $\Force^W(\set{\omega_n^W\st
n\in\omega})=\S{4.2}$.

\Corollary. If a model $W$ of set theory is absolutely definable in all forcing extensions (by the same formula, without parameters), then
$\Force^W=\S{4.2}$.

\Proof. If\/ $W$ is absolutely definable, then the $\omega_n^W$ are also absolutely definable, thereby avoiding in Theorem
\ref{Theorem.ValidWithParametersIsS4.2} the need for them to appear explicitly as parameters.\QED

An essentially identical argument shows, more generally, that if $W$ is absolutely definable from parameters in $X$, then $\Force^W(X)=\S{4.2}$. In
particular, no such $W$ is a model of the Maximality Principle. We now push the parameters a bit lower. Let $H_{\omega_2}$ denote the collection of
sets having hereditary size less than $\omega_2$.

\Theorem. For any model $W$ of set theory, $\Force^W(H_{\omega_2})=\S{4.2}$.\label{Theorem.Force(H_omega2)=S4.2}

\Proof. We will build an independent family of buttons and switches. In $W$, let $\omega_1=\bigsqcup_n S_n$ be a partition of $\omega_1$ into
infinitely many disjoint stationary subsets $S_n$. Let $b_n$ assert ``$S_n$ is not stationary''. Each assertion $b_n$ is false in $W$, since the $S_n$
are stationary there, but in any forcing extension, by shooting a club through the complement of any $S_n$, we can force $b_n$ to be necessary, while
preserving all stationary subsets of the complement of $S_n$ (see \cite[23.8, ex.~23.6]{Jech:SetTheory3rdEdition}). Thus, in any forcing extension of
$W$, we can push button $b_n$ without affecting any of the other buttons. (Note that this forcing collapses $\omega_1$ when $S_n$ is the sole remaining
stationary set on the list.) So the buttons $b_n$ are independent in $W$. For switches, let $s_k$ assert that the \GCH\ holds at $\aleph_{k+2}$. By
forcing over $W$ or any extension we can arrange the switches in any finite pattern, without adding subsets to $\omega_1^W$ and consequently without
affecting the buttons. So $W$ has an infinite independent family of buttons and switches using parameters in $H_{\omega_2}$, and so by Theorem
\ref{Theorem.IfButtonsThenS4.2} we conclude $\Force^W(H_{\omega_2}^W)=\S{4.2}$.\QED

This stationary set idea provides an alternative source of independent buttons for Lemma \ref{Lemma.ButtonsExist}, because if $\mathrm{V=L}$, then one
can use the $\mathrm{L}$-least partition of $\omega_1$ into $\omega$ many stationary sets; one advantage here is that this provides arbitrarily large
finite independent families of buttons that can be pushed without collapsing cardinals. Indeed, the button to collapse $\omega_1^\mathrm{L}$ is
equivalent to the infinite conjunction of these independent buttons.

Returning to Question \ref{Question.CanForceBeBetweenS4.2AndS5?}, we observe that if $W$ is a model of $\mathrm{V=L}$ or if $\Force^W=\S{4.2}$, then
clearly all the classes line up on the left side with $\Force^W(X)=\S{4.2}$. If $W\satisfies\MP$, then $\Force^W=\S5$ is on the right. If
$W\satisfies\MP(\R)$, then $\Force^W(\R)=\S5$. If $X$ has any uncountable parameter $x$, then the assertion $\psi=$ ``$x$ is countable'' is possibly
necessary but not true in $W$; consequently, $\possible\necessary\psi\implies\psi$ fails in $W$, and so $\Force^W(X)\neq\S5$. More generally, if $X$
has an element from which $\omega_1^W$ is absolutely definable, then $\Force^W(X)\neq \S5$. In fact, since \CH\ is always a switch and independent from
``$\omega_1^V$ is countable,'' we can conclude in this case that \axiomf{Dm} is not valid, so $\mathrm{Dm}\not\of\Force^W(X)$. If the set of parameters
supports long volume controls with independent switches, then $\Force^W(X)\of\S{4.3}$. If it supports many independent buttons and switches, then
$\Force^W(X)=\S{4.2}$. The situation $\S{4.2}\ofnoteq\Force^W(X)\ofnoteq\S5$ would occur if one could construct at least one button in $W$ using
parameters in $X$, but not arbitrarily large finite independent families of buttons and switches.

The hypothesis $\S5\of\Force^W(\R)$ is equivalent to $W\satisfies\MP(\R)$, which has large cardinal consistency strength. Specifically,
\cite{Hamkins2003:MaximalityPrinciple} shows it to be equiconsistent over \ZFC\ with the existence of a stationary proper class of inaccessible
cardinals. It is natural to inquire about the strength of weaker hypotheses concerning $\Force^W(\R)$. For example, we prove next that
$\mathrm{Dm}\of\Force^W(\R)$ already has large cardinal strength. Define that $\omega_1$ is {\df inaccessible to reals} if $\omega_1$ is an
inaccessible cardinal in $\mathrm{L}[x]$ for every real $x$ (equivalently, if $\omega_1^{\mathrm{L}[x]}<\omega_1$ for all $x\in\R$). The assertion that
there is a stationary proper class of inaccessible cardinals is expressed as a scheme in \ZFC, asserting that every class of ordinals (definable from
parameters) that is closed and unbounded in the class of all ordinals, contains an inaccessible cardinal.

\Theorem. If\/ $\mathrm{Dm}$ is valid for forcing with real parameters, then $\omega_1$ is inaccessible to reals and every $\mathrm{L}[x]$ for $x\in\R$
has a stationary proper class of inaccessible cardinals. Indeed, $\mathrm{L}_{\omega_1}\elesub \mathrm{L}$ and even $\mathrm{L}_{\omega_1}[x]\elesub
\mathrm{L}[x]$.\label{Theorem.IfDmInForceW(R)ThenLC}

\Proof. We show $\mathrm{L}_{\omega_1}\elesub \mathrm{L}$ by verifying the Tarski-Vaught criterion. This claim should be understood metatheoretically
and proved as a scheme. Suppose that $\mathrm{L}\satisfies\exists v\,\psi(u,v)$ where $u\in\mathrm{L}_{\omega_1}$. If there is no such $v$ inside
$\mathrm{L}_{\omega_1}$, then the assertion $\varphi(u)=$``there is $v\in\mathrm{L}_{\omega_1}$ such that $\mathrm{L}\satisfies\psi(u,v)$'' is false in
$W$. But it is certainly forceably necessary, because we could make $\varphi(u)$ true by collapsing cardinals until the least witness $v$ in
$\mathrm{L}$ is hereditarily countable. In other words, $\varphi(u)$ is a button in $W$. The parameter $u$, being in $\mathrm{L}_{\omega_1}$, is
hereditarily countable and can therefore be coded with a real. Meanwhile, the assertion \CH\ is a switch, and this switch is independent of
$\varphi(u)$, because the \CH\ and its negation can be forced over any model without collapsing $\omega_1$ and therefore without affecting the truth of
$\varphi(u)$. Since we have an independent button and switch in $W$ using a real parameter, it follows by Theorem
\ref{Theorem.ButtonAndSwitchCharacterizationsOfForce^W}, case (6), that \axiomf{Dm} is not valid in $W$ with real parameters, contrary to our
assumption that \axiomf{Dm} is valid for forcing over $W$ with real parameters.

By relativizing to any real $x$, we similarly conclude $\mathrm{L}_{\omega_1}[x]\elesub\mathrm{L}[x]$, again proved as a scheme. It follows by
elementary set theory that $\omega_1$ is (strongly) inaccessible in $\mathrm{L}[x]$. Also, if $C\of\ORD$ is a definable proper class club in
$\mathrm{L}[x]$, it follows that $C\intersect\omega_1$ is unbounded in $\omega_1$, and consequently $\omega_1\in C$. Thus, the class of inaccessible
cardinals in $\mathrm{L}[x]$ is a stationary proper class, as claimed.\QED

\Corollary. The following are equiconsistent over \ZFC:
\begin{enumerate}
 \item $\ZFC+\MP(\R)$.
 \item $\ZFC\,+$ \S5 with real parameters is valid for forcing.
 \item $\ZFC+\S{4W5}$ with real parameters is valid for forcing.
 \item $\ZFC+\theoryf{Dm.2}$ with real parameters is valid for forcing.
 \item $\ZFC+\theoryf{Dm}$ with real parameters is valid for forcing.
 \item $\ZFC\,+$ there is a stationary proper class of inaccessible cardinals.
\end{enumerate}\label{Corollary.EquiconsistenciesS5S4W5DmWithReals}

\Proof. The Maximality Principle $\MP(\R)$ asserts exactly that Axiom 5 is valid for forcing with real parameters, so $(1)\iff(2)$. Clearly
$(2)\implies(3)\implies(4)\implies(5)$, since $\S5\fo\S{4W5}\fo\theoryf{Dm.2}\fo\theoryf{Dm}$. In fact, $(4)\iff (5)$, since the axiom .2 is always
valid. Statement (5) implies, by Theorem \ref{Theorem.IfDmInForceW(R)ThenLC}, that (6) is true in $\mathrm{L}$. Finally, if there is a model where (6)
holds, then results in \cite{Hamkins2003:MaximalityPrinciple} produce a model satisfying $\ZFC+\MP(\R)$.\QED

The assertions in Corollary \ref{Corollary.EquiconsistenciesS5S4W5DmWithReals} should be understood as schemes asserting the relevant substitution
instances $\forall \vec x\in\R\ \varphi(\psi_0(\vec x),\ldots,\psi_n(\vec x))$.

\Question. Which modal theories arise as the valid principles of forcing with real parameters?

We know that $\S{4.2}\of\Force^W(\R)\of\S5$, and both of these endpoints can occur. But is there any model of set theory $W$ giving rise to any
intermediate modal theory?

\Section Restricting to a class of forcing notions \label{Section.RestrictingToAClass}

One can naturally restrict the scope of forcing extensions to those obtained by the members of a particular class $\Gamma$, so that
$\necessary_\Gamma\varphi$ means that $\varphi$ holds in all forcing extensions by forcing in $\Gamma$ and $\possible_\Gamma\varphi$ means that
$\varphi$ is forceable by some forcing in $\Gamma$. Natural classes $\Gamma$ would include \ccc\ forcing, proper forcing, and so on. The analogue of
the Main Question becomes:

\Question. For a given class of forcing $\Gamma$, what are the valid principles of\/ $\Gamma$ forcing? \label{Question.ValidPrinciplesForGamma}

There are many open questions here, which we leave for future projects. Let us close the paper by observing that even with \ccc\ forcing, the situation
changes dramatically.

\Theorem. If\/ \ZFC\ is consistent, then the \ZFC-provable principles of \ccc\ forcing do not include \S{4.2}.\label{Theorem.CCCisNotS4.2}

\Proof. It is easy to deduce in \S{4.2} the following Directedness axiom:
$$(\possible\necessary\varphi\wedge\possible\necessary\psi)\implies\possible\necessary(\varphi\wedge\psi).$$
We will show that the \ccc\ interpretation of this Directedness axiom is not valid in $\mathrm{L}$. This proof relies on the fact, a part of
mathematical folklore, that there are two fundamentally different ways to destroy a Suslin tree by \ccc\ forcing. If $T$ is a Suslin tree, then of
course forcing with $T$ itself adds a branch through $T$, making it non-Aronszajn. Alternatively, if $T$ is Suslin, then there is \ccc\ forcing making
$T$ into a special Aronszajn tree, a union of countably many antichains (see \cite[16.19 \& related]{Jech:SetTheory3rdEdition}). No further forcing can
add a branch through this tree without collapsing $\omega_1$, since the branch would have to contain at most one node from each antichain. Thus, these
two alternatives cannot be amalgamated by \ccc\ forcing. If $\varphi$ is the assertion ``the $\mathrm{L}$-least Suslin tree is not Aronszajn'' and
$\psi$ is the assertion ``the $\mathrm{L}$-least Suslin tree is special,'' then each of these statements is possibly necessary by \ccc\ forcing over
$\mathrm{L}$, but their conjunction is false in all \ccc\ extensions of $\mathrm{L}$. Thus,
$(\possibleccc\necessaryccc\varphi)\wedge(\possibleccc\necessaryccc\psi)$ holds in $\mathrm{L}$, but not
$\possibleccc\necessaryccc(\varphi\wedge\psi)$, violating Directedness.\QED

\Corollary. The same conclusion holds for any class $\Gamma$ of forcing notions containing all \ccc\ forcing, whose members in any \ccc\ extension
preserve $\omega_1$. This includes the classes of proper forcing, semi-proper forcing, cardinal-preserving forcing, and so on.

\Proof. In the proof of Theorem \ref{Theorem.CCCisNotS4.2}, the assertions $\varphi$ and $\psi$ are each \ccc-forceably necessary in $\mathrm{L}$ and
hence $\Gamma$-forceably necessary there, but the conjunction $\varphi\wedge\psi$ is not forceable by $\omega_1$-preserving forcing over $\mathrm{L}$,
and hence not $\Gamma$-forceable over $\mathrm{L}$. So the $\Gamma$-forcing validities in $\mathrm{L}$ do not include \S{4.2}.\QED

Nevertheless, one can easily verify that \S4 remains valid for \ccc\ forcing (and also for the other classes), and we conjecture that they do not go
beyond this.

\Conjecture. The \ZFC-provable principles of \ccc\ forcing are exactly \S4.\label{Conjecture.CCCvaliditiesAreS4}

A complete set of \S4 Kripke frames consists of finite pre-trees (partial pre-orders whose quotients are trees), and with them one might try to carry
out a similar analysis as in our Main Theorem, by finding set theoretical assertions to fulfill the Jankov-Fine assertions. The point is that branching
in these trees gives rise to behavior totally unlike either buttons or switches. Branching corresponds in set theory to the possibility of \ccc\
forcing extensions that cannot be amalgamated by further \ccc\ forcing, as in the folklore fact above, where one chooses either to specialize a Suslin
tree of $\mathrm{L}$ or to make it not Aronszajn. What is needed, therefore, is an elaborate generalization of this folklore idea, in which one can
successively make choices with \ccc\ forcing that cannot later be amalgamated by \ccc\ forcing, in such a way that every \ccc\ forcing extension is
included.

Under Martin's Axiom \MA, of course, the product of \ccc\ posets is again \ccc, and this implies that the Directedness Axiom is valid for \ccc\
forcing. As a consequence, all \S{4.2} assertions are valid for \ccc\ forcing over any model of \MA. We conjecture that this also is optimal. To prove
this, one would need to find a model of \MA\ with arbitrarily large finite independent families of \ccc\ buttons and switches.

\Conjecture. The $\ZFC+\MA$ provable principles of \ccc\ forcing are exactly \S{4.2}.

Lastly, we mention that for \ccc\ forcing, unlike the general situation with Corollary \ref{Corollary.EquiconsistenciesS5S4W5DmWithReals}, there is no
large cardinal strength to the hypothesis that \S5 is valid for \ccc\ forcing with real parameters. Specifically, \cite{Leibman2004:Dissertation}
proves that if \ZFC\ is consistent, then there is a model of \ZFC\ in which every \S5 assertion is valid for \ccc\ forcing with real parameters. If one
wants \S5 to be valid for \ccc\ forcing with parameters in $H_{2^\omega}$, however, then \cite{Hamkins2003:MaximalityPrinciple} shows that it is
equiconsistent, as in Corollary \ref{Corollary.EquiconsistenciesS5S4W5DmWithReals}, with a stationary proper class of inaccessible cardinals.

\bibliographystyle{alpha}
\bibliography{MathBiblio,HamkinsBiblio}

\end{document}

%% file: ArticleMacros.tex
\usepackage{latexsym,amsfonts,amsmath,amssymb}
\newtheorem{theorem}{Theorem}
\newtheorem{corollary}[theorem]{Corollary}
\newtheorem{sublemma}{Lemma}[theorem]
\newtheorem{lemma}[theorem]{Lemma}
\newtheorem{question}[theorem]{Question}
\newtheorem{observation}[theorem]{Observation}
\newtheorem{claim}[theorem]{Claim}
\newtheorem{subclaim}{Claim}[sublemma]
\newtheorem{conjecture}[theorem]{Conjecture}
\newtheorem{fact}[theorem]{Fact}
\newtheorem{definition}[theorem]{Definition}
\newtheorem{remark}[theorem]{Remark}
\newtheorem{example}[theorem]{Example}
\newtheorem{exercise}{Exercise}[section]
\def\Theorem #1.#2 #3\par{\def\claimname{Theorem}\setbox1=\hbox{#1}\ifdim\wd1=0pt
   \begin{theorem}{\rm #2} #3\end{theorem}\else
   \newtheorem{#1}[theorem]{#1}\begin{#1}\label{#1}{\rm #2} #3\end{#1}\fi}
\def\Corollary #1.#2 #3\par{\def\claimname{Corollary}\setbox1=\hbox{#1}\ifdim\wd1=0pt
   \begin{corollary}{\rm #2} #3\end{corollary}\else
   \newtheorem{#1}[theorem]{#1}\begin{#1}\label{#1}{\rm #2} #3\end{#1}\fi}
\def\Lemma #1.#2 #3\par{\def\claimname{Lemma}\setbox1=\hbox{#1}\ifdim\wd1=0pt
   \begin{lemma}{\rm #2} #3\end{lemma}\else
   \newtheorem{#1}[theorem]{#1}\begin{#1}\label{#1}{\rm #2} #3\end{#1}\fi}
\def\SubLemma #1.#2 #3\par{\def\claimname{Lemma}\setbox1=\hbox{#1}\ifdim\wd1=0pt
   \begin{sublemma}{\rm #2} #3\end{sublemma}\else
   \newtheorem{#1}{#1}[theorem]\begin{#1}\label{#1}{\rm #2} #3\end{#1}\fi}
\def\Question #1.#2 #3\par{\def\claimname{Question}\setbox1=\hbox{#1}\ifdim\wd1=0pt
   \begin{question}{\rm #2} #3\end{question}\else
   \newtheorem{#1}[theorem]{#1}\begin{#1}\label{#1}{\rm #2} #3\end{#1}\fi}
\def\Observation #1.#2 #3\par{\def\claimname{Observation}\setbox1=\hbox{#1}\ifdim\wd1=0pt
   \begin{observation}{\rm #2} #3\end{observation}\else
   \newtheorem{#1}[theorem]{#1}\begin{#1}\label{#1}{\rm #2} #3\end{#1}\fi}
\def\Claim #1.#2 #3\par{\def\claimname{Claim}\setbox1=\hbox{#1}\ifdim\wd1=0pt
   \begin{claim}{\rm #2} #3\end{claim}\else
   \newtheorem{#1}[theorem]{#1}\begin{#1}\label{#1}{\rm #2} #3\end{#1}\fi}
\def\SubClaim #1.#2 #3\par{\def\claimname{Claim}\setbox1=\hbox{#1}\ifdim\wd1=0pt
   \begin{subclaim}{\rm #2} #3\end{subclaim}\else
   \newtheorem{#1}{#1}[sublemma]\begin{#1}\label{#1}{\rm #2} #3\end{#1}\fi}
\def\Conjecture #1.#2 #3\par{\def\claimname{Conjecture}\setbox1=\hbox{#1}\ifdim\wd1=0pt
   \begin{conjecture}{\rm #2} #3\end{conjecture}\else
   \newtheorem{#1}[theorem]{#1}\begin{#1}\label{#1}{\rm #2} #3\end{#1}\fi}
\def\Fact #1.#2 #3\par{\def\claimname{Fact}\setbox1=\hbox{#1}\ifdim\wd1=0pt
   \begin{fact}{\rm #2} #3\end{fact}\else
   \newtheorem{#1}[theorem]{#1}\begin{#1}\label{#1}{\rm #2} #3\end{#1}\fi}
\def\Definition #1.#2 #3\par{\def\claimname{Definition}\setbox1=\hbox{#1}\ifdim\wd1=0pt
   \begin{definition}{\rm #2} {\rm #3}\end{definition}\else
   \newtheorem{#1}[theorem]{#1}\begin{#1}\label{#1}{\rm #2} {\rm #3}\end{#1}\fi}
\def\Remark #1.#2 #3\par{\def\claimname{Remark}\setbox1=\hbox{#1}\ifdim\wd1=0pt
   \begin{remark}{\rm #2} {\rm #3}\end{remark}\else
   \newtheorem{#1}[theorem]{#1}\begin{#1}\label{#1}{\rm #2} {\rm #3}\end{#1}\fi}
\def\Example #1.#2 #3\par{\def\claimname{Example}\setbox1=\hbox{#1}\ifdim\wd1=0pt
   \begin{example}{\rm #2} #3\end{example}\else
   \newtheorem{#1}[theorem]{#1}\begin{#1}\label{#1}{\rm #2} #3\end{#1}\fi}
\def\Exercise #1.#2 #3\par{\def\claimname{Exercise}\setbox1=\hbox{#1}\ifdim\wd1=0pt
   {\footnotesize\begin{exercise}{\rm #2} {\rm #3}\end{exercise}}\else
   \newtheorem{#1}[section]{#1}{\footnotesize\begin{#1}\label{#1}{\rm #2} {\rm #3}\end{#1}}\fi}
\def\QuietTheorem #1.#2 #3\par{\setbox1=\hbox{#1}\ifdim\wd1=0pt\proclaim{Theorem {\rm #2}}{#3}\else\proclaim{#1 {\rm #2}}{#3}\fi}
\newcommand{\proclaim}[2]{\smallskip\noindent{\bf #1} {\sl#2}\par\smallskip}
\def\Proclaim #1.#2 #3\par{\proclaim{#1 {\rm #2}}{#3}}
\newif\ifinproof\inprooffalse
\def\Proof#1. {\setbox1=\hbox{#1}\ifdim\wd1=0pt\begin{proof}\else\begin{proof}[#1]\fi\inprooftrue} % for use in amsart
\newcommand{\QED}{\end{proof}}

\def\BF#1.{{\bf #1.}}
\def\Abstract #1\par{\begin{abstract}#1\end{abstract}}
\def\Title #1\par{\title{#1}\maketitle}
\def\Author #1\par{\author{#1}}
\def\Acknowledgement#1\par{\thanks{#1}}
\def\Chapter #1\par{\chapter{#1}}
\def\Section #1\par{\section{#1}}
\def\QuietSection #1\par{\section*{#1}}
\def\SubSection #1\par{\subsection{#1}}
\def\SubSubSection #1\par{\subsubsection{#1}}
\def\MidTitle #1\par{\bigskip\goodbreak\centerline{\small\bf #1}\bigskip\noindent}
\def\Margin #1\par{\marginpar{\tiny #1}}

\newcommand{\citindex}{\index[citdep]}
\newcommand{\Ref}[2][Theorem]{\edef\entry{\ref{#2}..@#1 \ref{#2}!\ifinproof\arabic{chapter}.\arabic{section}.\arabic{theorem}..@\claimname\ \arabic{chapter}.\arabic{section}.\arabic{theorem}\else \arabic{chapter}.\arabic{section}..@Section \arabic{chapter}.\arabic{section}\fi|dotfill\ p.}\expandafter\citindex\expandafter{\entry}\ref{#2}}
\def\bottomnote #1\par{{\renewcommand{\thefootnote}{}\footnotetext{#1}}}
\newcommand{\F}{{\mathbb F}}

\renewcommand{\P}{{\mathbb P}}
\newcommand{\Q}{{\mathbb Q}}

\def\R{{\mathbb R}}

\newcommand{\Qdot}{{\dot\Q}}

\newcommand{\Rdot}{{\dot\R}}

\newfont{\msam}{msam10 at 12pt}

\newcommand{\of}{\subseteq}
\newcommand{\ofnoteq}{\subsetneq}
\newcommand{\fo}{\supseteq}

\newcommand{\set}[1]{\{\,{#1}\,\}}
\newcommand{\singleton}[1]{\left\{{#1}\right\}}

\newcommand{\elesub}{\prec}

\newcommand{\satisfies}{\models}
\newcommand{\forces}{\Vdash}
\newcommand{\proves}{\vdash}
\newcommand{\possible}{\mathop{\raisebox{-1pt}{$\Diamond$}}}
\newcommand{\necessary}{\mathop{\raisebox{-1pt}{$\Box$}}}

\newcommand{\necessaryccc}{\necessary_{\hbox{\scriptsize\ccc}}}
\newcommand{\possibleccc}{\possible_{\hbox{\scriptsize\ccc}}}

\newcommand{\cross}{\times}

\newcommand{\union}{\cup}

\newcommand{\intersect}{\cap}

\newcommand{\smalllt}{\mathrel{\mathchoice{\raise2pt\hbox{$\scriptstyle<$}}{\raise1pt\hbox{$\scriptstyle<$}}{\raise0pt\hbox{$\scriptscriptstyle<$}}{\scriptscriptstyle<}}}
\newcommand{\smallleq}{\mathrel{\mathchoice{\raise2pt\hbox{$\scriptstyle\leq$}}{\raise1pt\hbox{$\scriptstyle\leq$}}{\raise1pt\hbox{$\scriptscriptstyle\leq$}}{\scriptscriptstyle\leq}}}

\newcommand{\boolval}[1]{\mathopen{\lbrack\!\lbrack}\,#1\,\mathclose{\rbrack\!\rbrack}}
\def\[#1]{\boolval{#1}}

\newcommand{\UnderTilde}[1]{{\setbox1=\hbox{$#1$}\baselineskip=0pt\vtop{\hbox{$#1$}\hbox to\wd1{\hfil$\sim$\hfil}}}{}}
\newcommand{\Undertilde}[1]{{\setbox1=\hbox{$#1$}\baselineskip=0pt\vtop{\hbox{$#1$}\hbox to\wd1{\hfil$\scriptstyle\sim$\hfil}}}{}}
\newcommand{\undertilde}[1]{{\setbox1=\hbox{$#1$}\baselineskip=0pt\vtop{\hbox{$#1$}\hbox to\wd1{\hfil$\scriptscriptstyle\sim$\hfil}}}{}}
\newcommand{\UnderdTilde}[1]{{\setbox1=\hbox{$#1$}\baselineskip=0pt\vtop{\hbox{$#1$}\hbox to\wd1{\hfil$\approx$\hfil}}}{}}
\newcommand{\Underdtilde}[1]{{\setbox1=\hbox{$#1$}\baselineskip=0pt\vtop{\hbox{$#1$}\hbox to\wd1{\hfil\scriptsize$\approx$\hfil}}}{}}

\newcommand{\st}{\mid}
\renewcommand{\th}{{\hbox{\scriptsize th}}}

\newcommand{\minus}{\setminus}

\def\<#1>{\langle#1\rangle}

\newcommand{\ORD}{\mathop{\hbox{\sc ord}}}

\newcommand{\ZFC}{\hbox{\sc zfc}}

\newcommand{\CH}{\hbox{\sc ch}}

\newcommand{\GCH}{\hbox{\sc gch}}

\newcommand{\AC}{\hbox{\sc ac}}
\newcommand{\AD}{\hbox{\sc ad}}

\newcommand{\MA}{\hbox{\sc ma}}

\newcommand{\MP}{\hbox{\sc mp}}

\newcommand{\ccc}{{\hbox{\sc ccc}}}

\newcommand{\cell}[1]{\boxit{\hbox to 17pt{\strut\hfil$#1$\hfil}}}
\newcommand{\head}[2]{\lower2pt\vbox{\hbox{\strut\footnotesize\it\hskip3pt#2}\boxit{\cell#1}}}
\newcommand{\boxit}[1]{\setbox4=\hbox{\kern2pt#1\kern2pt}\hbox{\vrule\vbox{\hrule\kern2pt\box4\kern2pt\hrule}\vrule}}
\newcommand{\Col}[3]{\hbox{\vbox{\baselineskip=0pt\parskip=0pt\cell#1\cell#2\cell#3}}}
\newcommand{\tapenames}{\raise 5pt\vbox to .7in{\hbox to .8in{\it\hfill input: \strut}\vfill\hbox to
.8in{\it\hfill scratch: \strut}\vfill\hbox to .8in{\it\hfill output: \strut}}}
\newcommand{\Head}[4]{\lower2pt\vbox{\hbox to25pt{\strut\footnotesize\it\hfill#4\hfill}\boxit{\Col#1#2#3}}}
\newcommand{\Dots}{\raise 5pt\vbox to .7in{\hbox{\ $\cdots$\strut}\vfill\hbox{\ $\cdots$\strut}\vfill\hbox{\
$\cdots$\strut}}}
\renewcommand{\dots}{\raise5pt\hbox{\ $\cdots$}}
\newcommand{\factordiagramup}[6]{$$\begin{array}{ccc}
#1&\raise3pt\vbox{\hbox to60pt{\hfill$\scriptstyle
#2$\hfill}\vskip-6pt\hbox{$\vector(4,0){60}$}}&#3\\ \vbox
to30pt{}&\raise22pt\vtop{\hbox{$\vector(4,-3){60}$}\vskip-22pt\hbox
to60pt{\hfill$\scriptstyle #4\qquad$\hfill}}
     &\ \ \lower22pt\hbox{$\vector(0,3){45}$}\ {\scriptstyle #5}\\
\vbox to15pt{}&&#6\\
\end{array}$$}
\newcommand{\factordiagram}[6]{$$\begin{array}{ccc}
#1&&\\ \ \ \raise22pt\hbox{$\vector(0,-3){45}$}\ {\scriptstyle #2}
&\raise22pt\hbox{$\vector(2,-1){90}$}\raise5pt\llap{$\scriptstyle#3$\qquad\quad}&\vbox
to25pt{}\\ #4&\raise3pt\vbox{\hbox to90pt{\hfill$\scriptstyle
#5$\hfill}\vskip-6pt\hbox{$\vector(4,0){90}$}}&#6\\
\end{array}$$}
\newcommand{\df}{\it} % use italic for definition terms. Idea: also use this to create an index of definitions, if MakeIndex is true.